%% file: hhot181006.tex
\documentclass[12pt,a4paper]{article}
\usepackage{amsmath, amssymb, theorem, latexsym}
\usepackage{graphicx}

\newtheorem{theorem}{Theorem}[section]
\newtheorem{lemma}[theorem]{Lemma}
\newtheorem{proposition}[theorem]{Proposition}
\newtheorem{definition}[theorem]{Definition}
\newtheorem{assumption}[theorem]{Assumption}

\newtheorem{corollary}[theorem]{Corollary}

\newtheorem{remark}[theorem]{Remark}
\newtheorem{remarks}[theorem]{Remarks}

\allowdisplaybreaks
\def \eq {\ref}
\def\neweq#1{\begin{equation}\label{#1}}
\def\endeq{\end{equation}}
\newcommand\finbox{~\hfill$\Box$}%
\def \bbbr{\mathbb R}
\def\O{\Om}
\def\Om {{\Omega}}
\def\la {{\lambda}}

\def \Inte{{\rm Int\,}}

\def\pt#1{{\bf Proof of Theorem \ref{#1}.}}
\def\eq#1{(\ref{#1})}

\def\neweq{\begin{equation}}
\def\endeq{\end{equation}}

\def\endproof{{\mbox{}\nolinebreak\hfill\rule{2mm}{2mm}\par\medbreak} }
\newcommand {\ar}{\rightarrow}
\newcommand {\pa}{\partial}

\numberwithin{equation}{section}
\begin{document}
{\centering
\bfseries
{\Large Nodal domains and spectral minimal partitions.}
\\~\\~\\~\\

\par
\mdseries
\scshape
\small
B. Helffer$^{1,3}$\\
T. Hoffmann-Ostenhof$^{2,3}$ \\
S. Terracini$^4$
\par
\upshape
D\a'epartement de Math\a'ematiques, UMR 8628 Univ Paris-Sud$^1$\\
Institut f\"ur Theoretische Chemie, Universit\"at Wien$^2$\\
International Erwin Schr\"odinger Institute for Mathematical Physics$^3$\\
Universit\`a di Milano Bicocca, Dipartimento di Matematica e Applicazioni$^4$\\
\today

}

\begin{abstract}
We consider two-dimensional Schr\"odinger operators in bounded domains.
We analyze relations between the nodal domains, spectral minimal 
partitions  and spectral properties of the corresponding operator. The
main results concern the existence and regularity
 of the minimal partitions and the characterization
 of the minimal partitions associated with nodal sets
 as the nodal domains of Courant--sharp eigenfunctions.

\end{abstract}

\newpage 
\tableofcontents
\newpage

\section{Introduction and main results}\label{sintro}
We consider mainly two-dimensional Laplacians operators in bounded domains.
We would like to analyze the  relations between the nodal domains
 of the eigenfunctions of the Dirichlet Laplacians and the partitions by $k$
 open sets $D_i$ which are optimal in the sense that  the maximum over the $D_i$'s of the ground state energy of  the Dirichlet
 realization
 of the Laplacian in $D_i$ is minimal.\\
\subsection{Definitions and notations}
Let us consider a Schr\"odinger operator 
\begin{equation}\label{S}
H=-\Delta+V
\end{equation}
on a bounded  domain $\Om\subset \mathbb R^2$ with Dirichlet
boundary condition. 

In the whole article (except in Section \ref{susanna1}) , we will consider that $\Omega$ satisfies the
following
 condition of smoothness~:
\begin{assumption}\label{Asswhole1}~\\
$\Omega$ has compact and 
piecewise 
$\mathcal C^{1,+}$ boundary, i.e. piecewise $\mathcal C^{1,\alpha}$
for some $\alpha >0$. Moreover $\Omega$ satisfies the interior cone
property.
\end{assumption}
This allows a finite number of 
corners (and cracks) of opening $\alpha \pi$ (defined in Section \ref{spns}).\\
The other general assumption is that
\begin{assumption}\label{Asswhole2}~\\
The potential $V$ belongs to $ L^\infty(\Omega)$.
\end{assumption}

Under these assumptions (which will not be recalled at each
statement), 
$H$ is selfadjoint if viewed as the Friedrichs extension of the quadratic form associated to 
$H$ with form domain $W_0^{1,2}(\Om)$ and form core
$\mathcal C_0^\infty(\Om)$. We denote $H$ by $H(\Om)$. 
We are
interested in the eigenvalue problem for $H(\Om)$ and note that
under our assumptions $H(\Om)$ has compact
 resolvent
 and its spectrum, which will be denoted by  $\sigma(H(\Om))$ is
 discrete and consists of eigenvalues $\{\la_k\}_{k=1}^\infty$ 
with finite multiplicities which tend to infinity, so that 
\begin{equation}\label{lak}
\la_1<\la_2\le \la_3\le\dots\le \la_k\le  \dots
\end{equation}
and such that the associated eigenfunctions $u_k$ can be chosen to form an 
orthonormal basis for $L^2(\Om)$. \\
Without loss of generality we can assume that the $u_k$ are real 
valued and by elliptic regularity (see also Proposition \ref{holdercontinuous} in Section \ref{spns}) we have~:
\begin{equation}\label{regu}
u_k\in \mathcal C^{1,\alpha}(\Om)\cap \mathcal C_0^0(\overline\Om)\;,
\end{equation} 
for any $\alpha <1$.\\ 
We know that $u_1$ can be chosen to be strictly positive in $\Om$, but the other eigenfunctions 
$u_k$ ($k\geq 2$) must have zerosets. We define for any function $u\in
\mathcal C_0^0(\overline\Om)$
\begin{equation}
N(u)=\overline{\{x\in \Om\:\big|\: u(x)=0\}}
\end{equation}
and call the components of $\Om\setminus N(u)$ the nodal domains of $u$. The number of 
nodal domains of such a function will be called $\mu(u)$.

We now introduce the notions of partition and  minimal
 partition.

\begin{definition}\label{OPS}~\\
Let $1 \le k\in \mathbb N$. We
will call {\bf partition} (or $k$-partition if we want to indicate
 the cardinality of the partition) of $\Omega$ a family 
$\mathcal D=\{D_i\}_{i=1}^k$ of mutually disjoint subsets of $\Omega$~: 
\begin{equation}\label{DiO}
D_i\cap D_j=\emptyset\;,\; \forall i\neq j\mbox{\,and\,} 
 \cup_{i=1}^k D_i\subset \Om\;.
\end{equation}
We  call it {\bf open} if the $D_i$ are open sets of $\Omega$,
{\bf connected}
 if the $D_i$ are connected.\\
We denote by $\mathfrak O_k$ the set of open connected
$k$-partitions of $\Omega$.\end{definition}

We now introduce the notion of spectral minimal partition sequence.
\begin{definition}\label{regOm}~\\
Let $H=H(\Omega)$ as above. For $\mathcal D$ in $\mathfrak O_k$, we
introduce 
\begin{equation}\label{LaD}
\Lambda(\mathcal D)=\max_{i}\la(D_i)\;,
\end{equation}
where $\lambda(D_i)$ is the ground state energy of $H(D_j)$.
\end{definition}
\begin{remark}~\\
When $D$ is not sufficiently  regular, we define $\lambda(D)$
differently. See Definition~\ref{eigenvalue}.
\end{remark}

\begin{definition}~\\
For any integer $k\ge 1$, we define
\begin{equation}\label{frakL}
\mathfrak L_k=\inf_{\mathcal D\in \mathfrak O_k}\:\Lambda(\mathcal D).
\end{equation}
We call the  sequence $\{\mathfrak L_k\}_{k\ge 1}$ the {\em spectral minimal partition sequence} of 
$H(\Om)$. \\
For given $k$, we call $k$-- minimal 
partition a partition $\mathcal D\in \mathfrak O_k$,
such that  $\mathfrak L_k=\Lambda(\mathcal D)$. 
\end{definition}

\begin{remark}~\\
If $k=2$, it is rather well known (see for example \cite{HH:2005a} or
 \cite{CTV:2005}) that  $\mathfrak L_2$ is the second eigenvalue and
 the associated minimal  $2$-partition is the nodal partition
 associated to the second eigenfunction. 
\end{remark}

We now introduce the  notion of  strong partition.

\begin{definition}\label{def:reg_partition}~\\
A partition $\mathcal D=\{D_i\}_{i=1}^k$ of  $\Omega$ in $\mathfrak
O_k$ is called  \textbf{strong} if 
\begin{equation}\label{defstr}
\Inte(\overline{\cup_i D_i}) \setminus \pa \Om =\Om\;.
\end{equation}
\end{definition}

Attached to a  partition, we can naturally associate a closed
set in $\overline{\Omega}$ defined by
\begin{equation}\label{assclset}
N(\mathcal D) = \overline{ \cup_i \left( \partial D_i \cap \Omega
  \right)}\;.
\end{equation}

This leads us to introduce the definition  of  a  regular
 closed set. This definition is  
 is modeled on some (but not all) of  the properties of the nodal set
of an  eigenfunction of a Schr\"odinger operator (see Section \ref{spns}).

\begin{definition}\label{AMS}~\\
A closed set $N\subset \overline\Om$ is regular (and we write in this
case belongs to $N\in \mathcal M (\Om)$
(and we say in this case that it is regular) if $N$ meets the following requirements:
\\
\textbf{(i)}\\
There are finitely many distinct $x_i\in \Om\cap N$ and associated positive integers $\nu_i\text{ with }\nu_i\ge 2$
such that, in a sufficiently small neighborhood of each of the $x_i$, $N$ is the union of $\nu_i(x_i)$
$\mathcal C^{1,+}$ curves (non self-crossing) with one end at $x_i$ (and each pair 
defining at $x_i$ a positive angle in $(0,2\pi)$)
and such that in the complement of these points in $\Om$, $N$ is
locally diffeomorphic to a $\mathcal C^{1,1_-}$ (i.e. $\mathcal
C^{1,\alpha}$
 for any $\alpha \in ]0,1[$)  
curve.\\
\textbf{(ii)}\\
$\pa\Om\cap N$ consists of a (possibly empty) finite set of points
$z_i$, such that, at each $z_i$, $\rho_i$ $\mathcal C^{1,+}$ half-lines belonging to $N$ (with $\rho_i\geq
1$) hit the boundary. \\
\textbf{(iii)}\\
Moreover 
the half curves meet  with equal angle at each critical
 point of $N\cap \Omega $ and also at each point of $N\cap\partial  \Omega $
  together with the boundary.  
\end{definition}

Complementarily, we introduce the notion of regular partition.
\begin{definition}~\\
A strong partition $\mathcal D$ 
is regular (and we write in this case $\mathcal D \in \mathcal R(\Omega)$)
 if there exists a regular
closed
 set $N$ such that $\mathcal D =\mathcal D(N)$, where  $\mathcal D
 (N)$ is
 the family of the  connected components of $\Omega \setminus N$ belongs (by
 definition) to $\mathcal R(\Omega)$.
\end{definition}

 \subsection{Main results}
Although  some of the statements could be obtained under weaker
assumptions
 we assume below that $\Omega$ is bounded and connected.\\
 
It has been proved\footnote{But these papers
 treat only smoother boundaries than assumed in the whole article. So we will prove here a slight
 generalization.}
 by Conti-Terracini-Verzini \cite{CTV0, CTV2, CTV:2005} that 
\begin{theorem}\label{thstrreg}~\\
For any $k$, 
there exists a $k$-- minimal regular strong partition.
\end{theorem}

The first aim of this paper is to show the 
\begin{theorem}\label{thanyreg}~\\
Any $k$-- minimal partition has a connected, regular and strong representative.
\end{theorem}

Here we need
 to explain what we mean by representative (which involves implicitly
 the notion of capacity). This involves indeed a notion of equivalence 
 classes.
Two $k$-partitions $\mathcal D $ and $\widetilde{ \mathcal D}$
 are equivalent if there is a labeling such that
 for any ground state $u_i$ associated of $D_i$, there
 is a ground state $\widetilde u_i$ associated to $\widetilde D_i$ such that
$ u_i=\tilde u_i$ in $W_0^1(\Omega)$, and conversely.\\
Once this notion is introduced it is natural to look
 for the existence of a regular representative and uniqueness
 will always be inside this class.

In general, there is no reason for a minimal partition to be unique
 (and here we speak of uniqueness of equivalence classes). 
 This can for example occur in presence of symmetries. However, we
 will show that  the
 uniqueness property
 always holds for subpartitions of a given minimal  partition.
 More precisely, we have

\begin{theorem}\label{thsubpartition}~\\
Let $\mathcal D$ be a minimal $k$-partition relative to $\mathfrak
L_k(\Om)$.
 Let  $\mathcal D^\prime\subset \mathcal D$ be a  subpartition of
 $\mathcal D$ into $1\leq k^\prime< k$ elements   and  assume that 
$$\O':= \Inte\left(\bigcup_{D_i\in\mathcal D^\prime}\overline{D_i}\right)\;,$$
is connected.
Then  $\mathfrak L_k(\Om)=\mathfrak L_{k^\prime}(\O^\prime)$ and the
 $k'$-minimal partition of $\Omega'$ is unique.
\end{theorem}

A natural question is whether a minimal partition is the nodal partition induced by an eigenfunction.
Theorem \ref{partnod} gives a simple criterion for a partition to be 
 associated to a nodal set.  For this we need some additional
 definitions.

 We say that $D_i,D_j$ are \textbf{neighbors}
and we  write  $D_i\sim D_j$,  if the set \break $D_{i,j}:=\Inte(\overline {D_i\cup
  D_j})\setminus \pa \Om$ is connected. We then construct for  
each $\mathcal D$ a graph $G(\mathcal D)$ by associating to each $D_i$ a vertex and to each 
pair  of neighbors $(D_i,D_j)$ an edge. This is an undirected graph without
multiple edges or loops. We will say that the graph is \textbf{bipartite} if it
can be colored by two colors. We recall that the graph associated
 to a collection of nodal domains of an eigenfunction is always
 bipartite. In this case, we say that the partition is \textbf{admissible}. We have now the following converse theorem~:

\begin{theorem}\label{partnod}~\\
If the graph of the minimal partition of $\Omega$ 
 is bipartite,  this is a partition associated to the nodal set of an
 eigenfunction of $H(\Omega)$  corresponding to
 $\mathfrak L_k (\Omega)$.
\end{theorem}

This theorem was already obtained in \cite{HH:2005a} by adding a strong
a priori regularity assumption on the partition and the assumption that $\Omega$ is simply
connected. Any subpartition of cardinality  two $(D_i,D_j)$ 
 corresponds indeed to a
nodal partition of some eigenfunction associated 
 to the  second eigenvalue of $H(D_{i,j})$. This implies the Pair
 Compatibility Condition (see in Appendix \ref{sregcase})  and 
Theorem \ref{specconv} can be applied.\\
The proof given here is more general (but more difficult) and is
actually a byproduct of the proof of Theorem \ref{thanyreg}, which will
directly give
an eigenfunction whose nodal domains form the partition.

A natural question is now to determine how general is  the situation described in
the
 previous theorem. The surprise is that this will only occur in the so
 called Courant-sharp situation. Before stating
 precisely our secon main result we need to introduce some further statements  
 and notations. The  Courant Nodal
 Theorem says~: 
\begin{theorem}\label{Courant}~\\
Let $k\geq 1$, $\lambda_k=\lambda_k(\Omega)$ the $k$-th eigenvalue of $H(\Omega)$ 
 and $u$ any real associated  eigenfunction. Then the number of nodal domains $\mu(u)$ of $u$ satisfies 
$\mu (u)\le k$.
\end{theorem}
When the number of nodal domains $\mu(u)$ satisfies
 $$ \mu(u)=k\;,
$$
  we will say, as in 
\cite{AHH:2004}, that  $u$ is \textbf{Courant-sharp}.\\

\begin{definition}\label{LkGamma}~\\
For any integer $k\ge 1$, we denote by $L_k(\Omega)$
 (or simply $L_k$) the smallest eigenvalue whose eigenspace 
contains an eigenfunction with $k$ nodal domains.
\end{definition}
 In general, we will
 show in Corollary~\ref{rephr}, that 
\begin{equation}
\lambda_k (\Omega)\leq\mathfrak L_k(\Omega)\leq L_k (\Omega)\;.
\end{equation}
The last  goal consists in giving the full picture of the equality cases~:
\begin{theorem}\label{L=L}~\\
If $\mathfrak L_k(\Omega)=L_k(\Omega) $ or $\mathfrak L_k(\Omega)=\lambda_k(\Omega) $,  then 
$$\la_k(\Omega)=\mathfrak L_k(\Omega)=L_k (\Omega)\;.$$
 In addition, one can find in the eigenspace associated  to $\la_k$ an
 eigenfunction $u$ such that $\mu(u)=k$.
\end{theorem}
In other words, the only case when the $k$-- nodal domains
 of an eigenfunction of $H(\Omega)$ form a minimal partition
 is the case when this eigenfunction is Courant-sharp.

\subsection{Organization of the paper}
 The paper is organized as follows. We first start in Section
 \ref{spns}  by recalling and
 extending (at the boundary) results on the local
 properties of the nodal set of an eigenfunction. Section \ref{susanna1} is devoted
 to the analysis of the geometrical properties of minimal partitions in $\mathbb
 R^N$.
 Section \ref{susanna2} gives stronger results but limited to the
 two-dimensional case, which is our main subject. This gives in
 particular the proof of our first Main Theorem \ref{thanyreg}.
Sections \ref{sMNP} and \ref{spa} are devoted to additional properties of the
 minimal partitions. We discuss different notions related to the  spectrum
 and  revisit Pleijel's theorem
 and its proof. Section \ref{smpCs} gives the proof of the second Main
 Theorem \ref{L=L} permitting to show that when a $k$-minimal partition is
 a nodal family then the corresponding eigenvalue is the $k$-th
 one. In Section \ref{ssubp}, we complete the proofs and the
 statements concerning subpartitions. 
In Sections \ref{sdisk} and \ref{srectangle} we analyze in great
 detail the various spectra of
 specific $H(\Omega)$ in connection with minimal partitions. This leads in particular to nice 
 conjectures and open problems.
Finally, we develop in two appendices useful
results which will complete some proofs or  help
 the reader.\\

{\bf Acknowledgements}~\\
The two first authors are supported by the European Research Network
`Postdoctoral Training Program in Mathematical Analysis of Large
Quantum Systems' with contract number HPRN-CT-2002-00277, and the ESF
Scientific Programme in Spectral Theory and Partial Differential
Equations (SPECT). The support of the ESI was also very useful.\\
The two first authors   would also like to thank M. Van den Berg
 who was stimulating their interest for this question,
 M.~Dauge 
 for discussions on problems with corners 
 and V.~Bonnaillie-Noel and G.~Vial for their help in
 the analysis of the problem via numerical computations.
The third author was partially supported by  Italian MIUR, national project
``Variational Methods and Nonlinear Differential Equations''. She wishes to 
thank ESI, where a part of this work was done, for the support and the
kind hospitalty.

\section{Preliminaries: H\"older regularity of nodal sets}\label{spns}
It is a well known property of nodal sets of eigenfunctions to be the union of
curves ending either at interior singular points or at the boundary.  
This section is devoted to the analysis of the regularity of the nodal curves in the
  H\"older spaces  $\mathcal C^{1,\varepsilon}$. A word of caution must be entered at this point:
  with regularity we mean \emph {global regularity} of the nodal branch up to
  the singularities or the boundary. This is not a completely obvious issue  
 (basically because of the lack of regularity of our solutions and, possibly, of the boundary of the domain)
   and will require a revisitation of the       well known asymptotic estimates
  about critical points of eigenfunctions. To start with, we recall 
the classical local regularity result
 by Hartman and Wintner (\cite{hw1}, Corollary~1),  stating that  
interior critical points of non zero solutions to our class of equations
 are isolated and have finite (local) multiplicity $m$. In addition the 
 solution satisfies, for some $c\neq 0$, the asymptotic formula
\begin{equation}\label{asfo}
u(r,\theta)=c
r^{m+1}\cos((m+1)(\theta+\theta_0))+o(r^{m+1})\,,\qquad r=|z-z_0|\,.
\end{equation}
Here we identify $\mathbb R^2$ to $\mathbb C$ and use either
$z$, or $(x,y)$, or $(r,\theta)$ for a point of $\mathbb R^2$, with
the standard notations~:
$$
z= r \exp i \theta\;,\; z=x+i y\;,\; x= r\cos \theta, y=r\sin \theta\;.
$$

 We shall need a refined version of it which is stated below:

\begin{theorem}\label{reg_nodal_int}~\\
Let $\Omega$ be open and $V\in L^\infty(\Om)$. Assume $u\in W^{1,2}_{\rm loc}(\Omega)$ solves
$$ -\Delta u+V(x,y)u=0\;,$$ in the distributional sense. \\Let $z_0=(x_0,y_0)\in\Omega$ be such that $u(x_0,y_0)=0$ and $\nabla u(x_0,y_0)=0$; then, in a neighbourhood of $z_0$, 
\begin{itemize}
\item[(a)] There are an integer $n$, a  complex--valued function $\xi$
  of class $\mathcal C^{0,+}$ such that $\xi(z_0)\neq 0$ and
\begin{equation}
u_x+iu_y=r^n e^{-in\theta}\xi(x,y)\,,\qquad r=|z-z_0|
\end{equation}
\item[(b)] There is a  function $\tilde \xi$ of class $\mathcal C^{0,+}$ such that $\tilde \xi(x_0,y_0)=0$ and
\begin{equation}
u(x,y)=\dfrac{r^{n+1}}{n+1}\left(\Re(\xi(x_0,y_0))\cos{(n+1)\theta}+
\Im(\xi(x_0,y_0))\sin{(n+1)\theta}+\tilde\xi(x,y)\right).
\end{equation}
\item[(c)]  There exists a positive radius $R$ such that $u^{-1}(\{0\})\cap B(z_0,R)$
 is composed by $2n$ $\mathcal C^{1,+}$-simple arcs which all end in $z_0$ 
and whose tangent lines at $z_0$ divide the disc into $2n$ angles of equal amplitude. 
\end{itemize}
\end{theorem}

\begin{proof}~\\
Following the paper by Hartman and Wintner \cite{hw1}, denote $w=u_y+iu_x$ and set $z_0=0$. Then, it is shown there that, if,
\begin{equation}\label{littleo}
u=o(|z|^{k}),
\end{equation}
for some integer $k\geq 0$, then the Cauchy formula is available:
\begin{equation}\label{cauchy}
2\pi i\frac{w(\zeta)}{\zeta^k}=\int_{|z|=R} \frac{w(z)}{z^k(z-\zeta)}dz-\int_{|z|<R} \frac{V(z)u(z)}{z^k(z-\zeta)}dx\,dy.
\end{equation}
where $R>0$ is fixed and the double integral over the disk is
 absolutely convergent. We now wish to show that it follows that the
 left-hand side is
 H\"older continuous in $\zeta$ in a neighbourhood of the origin. 
The line integral  is smooth in $\zeta$, since the integrand has no 
singularities on the circle. Concerning the second term, notice that 
we can find a constant $K$ such that
\begin{equation*}
\begin{split}
\left|\int_{|z|<R} \frac{V(z)u(z)}{z^k} \left(\frac{1}{z-\zeta_1}-\frac{1}{z-\zeta_2}\right)dx\,dy\right|\leq
\\
\int_{|z|<R} \frac{|V(z)u(z)|}{|z|^k} \left|\frac{|\zeta_1-\zeta_2|}{|z-\zeta_1| |z-\zeta_2|}\right|dx\,dy\leq \\
K |\zeta_1-\zeta_2| \vert\log|\zeta_1-\zeta_2|\vert.
\end{split}
\end{equation*}
Now we wish to show that \eq{littleo} can not be verified for every
integer. To this aim, we integrate  Equation \eq{cauchy} over the disk and, taking  absolute values, we obtain: 
\begin{equation}\label{cauchyint}
2\pi \int_{|z|<R}\frac{|w(z)|}{|z|^k} dx\, dy \;\leq 2\pi
R\int_{|z|=R} \frac{|w(z)|}{|z|^k}|dz|+2\pi  R\int_{|z|<R} \frac{|V(z)| |u(z)|}{|z|^k}dx\,dy.
\end{equation}
Following \cite{hw1}, and using the identity
\begin{equation}\label{identi}
u(r,\theta)=\int_0^r \left (u_x(\rho, \theta)\cos\theta+u_y(\rho, \theta)\sin\theta\right)d\rho\,,
\end{equation}
 we observe that~:
\[
|u(z)|\leq \int_0^1|z w(tz)|\, dt\;,
\]
implies 
\[
\int_{|z|<R} \frac{|V(z)| |u(z)|}{|z|^k}dx\,dy\leq K\int_{|z|<R} \frac{ |z w(z)|}{|z|^k}dx\,dy\leq KR\int_{|z|<R} \frac{ |w(z)|}{|z|^k}dx\,dy\;.
\]
Thus, for $R$ sufficiently small, inequality \eqref{cauchyint} leads to
\begin{equation}\label{cauchyint2}
\int_{|z|<R}\frac{|w(z)|}{|z|^k}dx\,dy \leq2R\int_{|z|=R} \frac{|w(z)|}{|z|^k}|dz|.
\end{equation}
We have now fixed $R>0$ such that \eqref{cauchyint2} is satisfied. Let us assume that $w(z_0)\neq 0$ for some $|z_0|<R$. Then, for a constant independent of $k$ there holds
\[|w(z_0)| \leq K \left(\frac{|z_0|}{R}\right)^k\;,\qquad k=1,2,\dots\;.
\]
Let us take the limit $k\rightarrow +\infty$ in this inequality. Then the limit of the sequence $(|z_0|/R)^k$ does not vanish, in contradiction with $|z_0|<R$. This completes the proof of point $(a)$ in the statement of the Theorem.
Point $(b)$ follows from point $(a)$ together with  the identity \eqref{identi}.

To prove point $(c)$ we choose a branch of the nodal set and we choose, as a regular parametrization the path $z(t)=r(t)e^{i\theta(t)}$, where the pair $(r(t),\theta(t))$ solves the following system of ordinary differential equations:
\begin{equation*}
\begin{cases}
\dot r=\dfrac{1}{r^{n+1}}\left(xu_y-yu_x\right) & \\
\dot \theta=\dfrac{1}{r^{n+2}}\left(xu_x+yu_y\right)\;.,& 
\end{cases}
\end{equation*}
where $\dot r$ and $\dot \theta$ denote respectively the derivative of
$r(t)$ and $\theta(t)$ with respect to $t$.\\
One can easily prove by the points $(a)$ and $(b)$ that both functions 
$t\mapsto \dot r(t)$ and $t\mapsto r(t)\dot \theta (t)$ are H\"older continuous; therefore both
$r$ and $\theta$ are H\"older continuous functions. Hence they can be extended through 
the singularity. Since the parametrization is regular ($\dot z\neq 0$), the assertion follows from the equation
$$ \dot z= \dot r(t)e^{i\theta(t)}+ir(t)\dot
\theta(t)e^{i\theta(t)}\;.
$$
\end{proof}

\begin{remark}\label{singular_potential}~\\
 Theorem \ref{reg_nodal_int} extends, with the same argument, to the case when the potential $V$ has a singularity at $z_0$, provided there exists $\beta<1$ and $K$ such that
\[|V(x,y)|\leq  \dfrac{K}{|z-z_0|^\beta}.\]
This fact will be useful when we shall consider the case of domains
with corners or cracks; indeed such singular potentials result as
conformal factors associated with the complex exponentials. \\
Note that this singular situation was also analyzed, but for the
interior problem,  in \cite{HO2} and \cite{HO2NA} and that
 in this case the authors obtain a better regularity.
\end{remark}

In order to examine the regularity up to the boundary of the nodal partition associated to an eigenfunction we now  extend a known result by Alessandrini \cite{al:nodal} (which treats the
convex case)  to our setting. The proof exploits the classical Kellog--Warschawski theorem on the boundary regularity of conformal mappings which states that any conformal map on a $\mathcal C^{1,\varepsilon}$ domain extends continuously on the boundary keeping the same regularity
(see the book by Pommerenke \cite{Pom}, Theorem 3.6 in Chapter~3).

\begin{theorem}\label{reg_nodal}~\\
Let $\varepsilon >0$ and $\Omega$ be an open set with $\mathcal C^{1,\varepsilon}$ boundary and $V\in L^\infty(\Om)$. Assume $u\in W^{1,2}_0(\Omega)$ solves
$$ -\Delta u+V(x,y)u=0\;,$$ in the distributional sense. Then the associated nodal partition 
is regular. More precisely if $\overline{u^{-1}(\{0\})}$
 intersects $\partial\Omega$ at $z_0$, 
then there exist an integer $m$ and $R>0$ such that
 $\overline{u^{-1}(\{0\})}\cap B(z_0,R)$
 is composed by $m$ $\mathcal  C^{1,\varepsilon}$-simple arcs which all end in $z_0$ 
and whose tangent lines at $z_0$ divide the tangent cone $\Gamma(z_0)$
 into $m+1$ angles of equal opening. 
\end{theorem}
\begin{proof}~\\
The result immediately follows from Theorem \ref{reg_nodal_int} in the case of the half--plane: indeed one can extend $u$ by a reflection to the other half--plane and reduce to the case of the interior zeros. The general case reduces to that of the half--space through the Riemann  mapping Theorem. Indeed,  by \cite{Pom} (Theorem 3.6 in Chapter 3) the H\"older regularity $\mathcal  C^{1,\varepsilon}$ of $\partial \Omega$ implies the same regularity property for the extensions, up to boundary, of the Riemann map $f$ and of its inverse. Since the composition of $\mathcal  C^{1,\varepsilon}$ maps enjoys the same regularity property, the statement follows.
\end{proof}
Now we wish to extend Theorem \ref{reg_nodal} to the case of domains possessing corners or cracks. To be precise we start with the following

\begin{definition}~\\ Let $\varepsilon \in ]0,1]$. We say that $\partial\Omega$ has a $\mathcal C^{1,\varepsilon}$--corner of opening $\alpha\pi$ ($0\leq \alpha\leq 2$) at $z_0$ if, in a sufficiently small neighborhood, $\partial\Omega$ contains the union of two
curves of class $\mathcal C^{1,\varepsilon}$ (non self-crossing)
ending at $z_0$, and such that $\Omega$ lies in the curvilinear sector
of angle opening $\alpha\pi$ spanned by the two arcs, which does not intersect other components of the boundary $\partial\Omega$. 
\end{definition}

\begin{remark}~\\
Note that the boundary $\partial\Omega$ can have several corners of 
angle opening $\alpha_i$ at 
the same point $z_0$: of course the sum of all the angles
 does not exceed $2\pi$.  Moreover, 
it is worthwile noticing that we allow the presence of cracks
(i.e. corners of 
 angle opening $2\pi$ where the two curves coincide), exterior cusps
 (i.e. corners of angle opening $2\pi$ spanned by two distinct
 curves), as well as angles of any possible positive angles,
 (positiveness is required by the interior cone property). Finally, a
 corner can be a point of smoothness of the boundary, when its angle opening is $\pi$.
\end{remark}

Our next goal is to prove the following result

\begin{theorem}\label{reg_nodal_corner}~\\
Let $V\in L^\infty(\Om)$ and $\varepsilon \in ]0,1]$. Assume $u \in W_0^{1,2}(\Omega)$ solves
$-\Delta u+V(x,y)u=0$ in the distributional sense, in a neigbourhood of some $z_0\in\partial\Omega$, a $\mathcal C^{1,\varepsilon}$--corner of opening $\alpha\pi$ ($0< \alpha\leq 2$). If $\overline{u^{-1}(\{0\})}$
 intersects $\partial\Omega$ at $z_0$, 
then there exist an integer $m$ and $R>0$ such that
 $\overline{u^{-1}(\{0\})}\cap B(z_0,R)$
 is composed by $m$ $\mathcal C^{1,\varepsilon'(\alpha)}$-simple arcs which all end at $z_0$ 
and whose tangent lines at $z_0$ divide the tangent cone $\Gamma(z_0)$
 into $m+1$ angles of equal amplitude. In addition
\[\varepsilon'(\alpha)=\begin{cases}\varepsilon\min(\alpha,1/\alpha)\qquad &\mbox{if } 1/2<\alpha\leq 2 \\
2^n \varepsilon\alpha\qquad &\mbox{if } 1/2^{(n+1)}<\alpha\leq 1/2^n
\end{cases}
\]\end{theorem}

To prove the theorem we shall first straighten the corner and then apply Theorem \ref{reg_nodal}. We shall need the following basic result.

\begin{proposition}\label{straighten}~\\
 Let $\varepsilon \in ]0,1]$ and let $C$ be a H\"older--continuous arc ending at the origin, without self--intersections. Let $w(\tau)$, $\tau\in[0,\overline \tau]$ be a regular parametrization of $C$ such that $w(0)=0$ and $w'(0)\neq 0$, and define
the curve $\mathcal C^{1/\alpha}$  by the parametrization
\[t \mapsto v(t):=\left(w(t^\alpha)\right)^{1/\alpha}.\]
Then, for any $\alpha>0$
\[C\in \mathcal C^{1,\varepsilon}\quad\Longrightarrow\quad C^{1/\alpha}\in \mathcal C^{1,\varepsilon\min(1,\alpha)}\]
\end{proposition}
\begin{proof}~\\
We have 
\[v'(t)=w'(t^\alpha)u(t^\alpha),\qquad u(\tau):= \left(\tau w(\tau)\right)^{-1+1/\alpha}. 
\]
Obviously $v$ defines a regular parametrization of $C^{1/\alpha}$. At first we remark that $u$ is H\"older continuous with exponent $\varepsilon$. Indeed
\[\left|\dfrac{w(\tau_1)}{\tau_1}-\dfrac{w(\tau_2)}{\tau_2}\right|=\left|\int_0^1(w'(\tau_2
s)-w'(\tau_1 s))\,ds\right|\leq K |\tau_1-\tau_2|^\varepsilon.\]
Therefore the product $w'(\tau)u(\tau)$ is of the same class $\mathcal C^{0,\varepsilon}$ and the composition 
$v'(t)=w'(t^\alpha)u(t^\alpha)$ is in the H\"older space $\mathcal C^{0,\varepsilon\min(1,\alpha)}$.
\end{proof}

\begin{proof}~\\
We are in  position to prove Theorem \ref{reg_nodal_corner}. We 
 consider separately the two cases~:\\

 First, we assume the case of openings satisfying the inequality $1/2<\alpha\leq 2$. We  straighten the corner as in Proposition \ref{straighten}, by the map $z\to z^{1/\alpha}$. Then, through this composition, the boundary $\partial \Omega^{1/\alpha}$  becomes smooth (of class $\mathcal C^{1,\varepsilon\min(1,\alpha)}$) while the potential $V$ has to be multiplied by the conformal factor  $2\alpha^2|z-z_0|^{2(\alpha-1)}$, which is singular whenever $\alpha<1$. As already
observed in Remark \ref{singular_potential}, this is not a problem if
$\alpha>1/2$. Thanks to Theorem \ref{reg_nodal}, the nodal set
 of the composition $z\mapsto u(z^\alpha)$ is the union of arcs of
 class $\mathcal  C^{1,\varepsilon\min(1,\alpha)}$. Now we take its
 inverse image through the map $z\to z^{\alpha}$ and, applying again
 Proposition \ref{straighten} we obtain the desired value of
 $\varepsilon'(\alpha)$.\\

 Next we turn to the case when the opening is
 too small, that is when $1/2^{(n+1)}<\alpha\leq 1/2^n$, for some
 $n\geq 1$. 
Using again  Theorem 3.6 in Chapter 3 of \cite{Pom}, one can easily
 construct, locally in $\Omega$, a conformal map of class $\mathcal
 C^{1,\varepsilon}$ up to the boundary such that the image of one of
 the two arcs is a straight segment. Next step is to reflect the
 domain about this line and extend the function on  the reflected
 corner, in such a way to double the opening, which is now
 $2\alpha$. In this procedure, the second arc, being  composed with a
 $\mathcal C^{1,\varepsilon}$ map, still remains in the same H\"older class. We iterate this reflection procedure $n$ times, until $1/2<2^n\alpha\leq 1$, and we afterward proceed as in the proof of the case $\alpha >1/2$.  
\end{proof}

Using the same technique of straightening the angles by conformal maps, one can easily prove the following
\begin{proposition}\label{holdercontinuous}~\\
Let $V\in L^\infty(\Om)$ and $\varepsilon \in ]0,1]$. Assume $u\in W^{1,2}_0$ solves
$-\Delta u+V(x,y)u=0$ in the distributional sense, in a neigbourhood of some 
$z_0\in\partial\Omega$, a $\mathcal C^{1,\varepsilon}$--corner of opening 
$\alpha\pi$ ($0< \alpha\leq 2$). 
Then, if $\alpha \leq 1$, $u\in\mathcal C^{1,\varepsilon}$, locally at $z_0$;  
otherwise, if $1<\alpha\leq 2$,  we only have $u\in\mathcal C^{0,1/\alpha}$.
\end{proposition}

\newpage
\input susanna06m.tex

\newpage
\section{More  on nodal sets and partitions}\label{sMNP}

We continue by discussing more deeply  the links between the various
spectral sequences.\\

The first  important property is given by~:
\begin{proposition}\label{l<L}~\\
Let $H(\Om)$ be defined as above. Then 
\begin{equation}\label{frakL<}
\mathfrak L_k (\Omega)<\mathfrak L_{k+1}(\Omega)\text{ for } k\ge 1.
\end{equation}
\end{proposition}
{\bf Proof}~\\
We take indeed a minimal $(k+1)$-partition of $\Omega$. We have
 proved that this partition is regular. If we take any
 subpartition by $k$ elements of the previous partitions. This
 cannot be a $k$-minimal partition (it has not the ``strong
 partition'' property). So the inequality in \eqref{frakL<} should be strict.

The second property concerns the {\bf domain monotonicity}~\\
It is indeed immediate to verify~:
\begin{proposition}~\\
If $\Omega \subset \widetilde \Omega$, then
$$ \mathfrak L_k (\widetilde \Omega) \leq \mathfrak L_k(\Omega)\;,\; \forall
k\geq 1\,.$$
\end{proposition}
We observe indeed that each  partition of  $\Omega$ 
 is a  partition of $\widetilde \Omega$.

\begin{remark}~\\
The analysis of the equality in the proposition will involve the
capacity
 of $\widetilde \Omega \setminus \Omega$. See \cite{AHH:2004}
 and references therein.
\end{remark}

We now come back to  a definition which was briefly
mentioned in the introduction. Having in mind
 Definition \ref{LkGamma},  we denote, for any integer $k\ge 1$, by $L_k(\Omega) $ the smallest eigenvalue whose eigenspace 
contains an eigenfunction with $k$ nodal domains. We take 
$L_k=+\infty$,  if there are no 
eigenfunctions with $k$ nodal domains. We call this sequence the spectral nodal sequence.

\begin{proposition}\label{LkmathcalLk}~\\
Let $\lambda$ be an eigenvalue corresponding to an eigenfunction  with  $k$ nodal domains.
Then
\begin{equation}\label{ineq>}
 \mathfrak L_k  \leq \lambda\;.
\end{equation}

\end{proposition}

{\bf Proof of Proposition \ref{LkmathcalLk}}.\\
If $u$ is an eigenfunction associated with $\lambda$ and with
 $k$ nodal domains, then, taking as $\mathcal B_0$ the collection of these nodal domains, we obtain~:
\begin{equation}
\underset{\mathcal B\in\mathfrak O_k}\inf \Lambda(\mathcal B)\leq \Lambda(\mathcal B_0)  \leq \lambda\;.
\end{equation}

\begin{proposition}\label{lambdakLk}~\\
Let $\lambda =\lambda_k$ be an eigenvalue  $H(\Omega)$. 
Then
\begin{equation}\label{ineq<}
\lambda_k  \leq \underset{\mathcal B\in\mathfrak O_k}\inf \Lambda(\mathcal B).
\end{equation}
\end{proposition}

{\bf Proof of Proposition \ref{lambdakLk}}.\\
The basic idea (which is already present in Courant's Theorem) is simply the following.
 We can assume (using Proposition \ref{l<L}) 
 $$ \lambda_{k-1} <\lambda_k\;.$$
Attached to a minimal  $\mathcal B_k$ (hence regular), 
we have a $k$-dimensional space
 in $H_0^1(\Omega)$ generated by the ground states
 of the $D_i$ ($i=1,\dots, k$). We can find in this space a non
 trivial
 element which is orthogonal to the eigenspace
 corresponding to the eigenvalues which are $\leq \lambda_{k-1}$,
 whose energy is $\mathfrak L_k$, hence by the Minimax Principle
$\lambda_k \leq \mathfrak L_k$.\\
Suppose now that we have the equality $\lambda_k = \mathfrak L_k$.
Again by the proof of the Minimax Principle, this non trivial element
should be an eigenfunction which is consequently Courant-sharp
 and we have consequently $\lambda_k=L_k=\mathfrak L_k$.\\

The following corollary  is just a rephrasing of Propositions 
\ref{lambdakLk} and \ref{LkmathcalLk}.
\begin{corollary}\label{rephr}~\\
We have
\begin{equation}\label{ineqbetw}
L_k\ge \mathfrak L_k\ge \la_k\;,\, \forall k \geq 1\;.
\end{equation}
\end{corollary}

In particular, if $L_k=\lambda_k$ (also called the Courant-sharp case
 see in \cite{AHH:2004})
 the nodal domain of a corresponding eigenfunction gives a minimal partition.

\begin{remarks}\label{S(H)}~
\begin{enumerate}
\item
For the one dimensional case the
standard Sturm-Liouville theory leads easily to the following 
\begin{equation}
L_k = \mathfrak  L_k = \lambda_k\;,\; \forall k \geq 1\;.
\end{equation}
\item
It is easy to show,  that for a given $H$ 
\begin{equation}\label{Ssigma}
\mathfrak L_1=L_1=\lambda_1,
\end{equation}
(by the property of the ground state)
and that 
\begin{equation}
 \: \mathfrak L_2=L_2=\lambda_2,
\end{equation}
by the orthogonality of $u_2$ to the ground state combined with
Courant's nodal Theorem.  (See also \cite{CTV:2005}, Corollary 4.1  (Case
$V=0$)), 
but the extension to $V\in L^\infty$ is not  a problem.
\item
The sequence $L_k$ is not necessarily monotone : see for example \eqref{L3lambda}.
\end{enumerate}
\end{remarks}

One also observes
 that,  using  (\ref{ineqbetw}) and  the property that 
 $\lambda_k \ar +\infty$,
\begin{equation}
\lim_{k\ar +\infty} \mathfrak L_k = +\infty\;.
\end{equation}

\newpage 
\section{Playing around Pleijel's argument.}\label{spa}

It is a well known result of Pleijel that the one dimensional result that
 the $k$-th eigenvector of a Sturm-Liouville operator on an interval
 has only $k$ nodal domains 
 cannot be extended to higher dimension. The $k$-th Courant-sharp eigenfunctions
 (i.e. $k$-th eigenfunctions with $k$ nodal domains) can only be found for a finite number of $k$'s.\\
We will show in this section, that the arguments behind the proof of this theorem give also many informations on the spectral minimal partition sequence
 in comparison with the spectral sequence and the nodal sequence.\\

Let us look at 
a universal  lower bound for $\mathfrak L_k (\Omega)$. We actually obtain~:
\begin{proposition}~\\
Considering the Dirichlet Laplacian, we have
$$
\mathfrak L_k(\Omega) \geq k  \frac{\pi j^2}{|\Omega|}\;,
$$
where $|\Omega|$ denotes the area of $\Omega$ and $j$ is the smallest positive $0$ of the Bessel function $J_0$~:
\begin{equation}\label{valueofj}
j\sim 2. 4048\dots\;.
\end{equation}
\end{proposition}

{\bf Proof}~\\
For any $D_j$ of a partition, we have by the Faber-Krahn inequality 
$$
|D_j| \lambda(D_j) \geq \pi j^2\;.
$$
The Faber-Krahn inequality gives indeed ~:
\begin{equation}
\lambda(D) \geq \frac{\lambda(B_{1/\sqrt{\pi}})}{|D|}\;,
\end{equation}
for any open set $D$.\\
The lowest eigenvalue for the disk of radius $1$ is known to be :
\begin{equation}
\lambda(B_1)=   j^2\;,\;\mbox{with\;} \pi j^2 \sim 18.1695\;.
\end{equation}

Summing up over $j$, we obtain
\begin{equation} \label{upb}
\pi j^2 k \leq \sum_j |D_j| \lambda(D_j) \leq |\Omega| \max \lambda(D_j)\;.
\end{equation}
Taking the infimum over the partition leads to the result.

\begin{remark}~\\
Using Corollary \ref{rephr}, this implies
$$
L_k(\Omega) \geq k  \frac{\pi j^2}{|\Omega|}\;.
$$
\end{remark}

We conclude this section with a classical result of Pleijel, \cite{Pleijel:1956}:
\begin{theorem}\label{Pleijel}~\\
If $\Omega$ is bounded and  smooth, 
\begin{equation}\label{Plej}
\{k\;|\; u_k\text{ has }k\text{ nodal domains}\} \text{ is finite.}
\end{equation}
\end{theorem}

This holds in larger generality for bounded potentials and also for 
higher dimensions. 

Let us describe for completeness how  Pleijel's Theorem is proved. The Weyl theory says that
\begin{equation}\label{Wey}
\lambda_n \sim \frac{4 \pi n}{|\Omega|}\;,
\end{equation}
as $n\ar + \infty$.

If $u_n$ is an eigenfunction associated to $\lambda_n$ with $n$ nodal
 domains, 
 we obtain immediately a contradiction for $n$ large between
 (\ref{Wey}) and (\ref{upb}) (applied with the family of nodal domains
 of $u_n$),
 having in mind the value of $j$ given in (\ref{valueofj}). So $u_n$
 cannot have $n$ nodal domains ! More precisely, if there exists a smallest
  $n(k)$ such that $\lambda_{n(k)} = L_k$, we obtain asymptotically
\begin{equation}
\liminf_{k\ar +\infty} \frac{n(k)}{k} \geq j^2/4 >1\;.
\end{equation}

A more difficult  question is to determine whether  $L_k$ is always finite.

\begin{proposition}\label{propcomphexa}~\\
In the case of the Laplacian and if $\Omega$ is regular, we have\footnote{Thanks to M. Van den
  Berg for discussions. In particular he conjectures the existence
 of the limit $\lim_{k\ar +\infty} \frac{\mathfrak L_k}{k}$ and that the limit
 is actually $\lambda_1 (\mbox{Hx}_1)$.} 
\begin{equation}\label{comphexa}
\limsup_{k\ar +\infty}  \frac{\mathfrak L_k(\Omega)}{k} \leq 
\lambda_1 (\mbox{Hx}_1) / |\Omega|\;,
\end{equation}
where $\mbox{Hx}_1$ is the regular hexagon of area $1$.
\end{proposition}
For the proof, we just use a (non strong) partition  of $\Omega$ by equal hexagons
of area at most $\frac{|\Omega|}{ k}$.

\begin{remark}.\\
Adding a potential does not create any diffuculty and the previous discussions 
can be easilily adapted to go from $H_0(\Omega)$ to $H(\Omega)$.
Concerning the values $L_k$ we obtain immediately,
\begin{equation}
L_k^V(\Omega) \geq k \frac{\pi j^2}{|\Omega|}  - \sup |V|\;.
\end{equation}

On the other hand,  using the minimax, there are no problem to show that
\begin{equation}
\mathfrak L_k^V(\Omega) \geq k \frac{\pi j^2}{|\Omega|}  - \sup |V|\;.
\end{equation}
This implies in particular that
\begin{equation}
\liminf_{k\ar +\infty}  \frac{\mathfrak L_k^V(\Omega)}{k} \geq 
\frac{\lambda_1 (B_{1/\sqrt{\pi}}) }{ |\Omega|}> \frac{4\pi}{|\Omega|}\;,
\end{equation}
is satisfied in full generality.
\end{remark}

\newpage 
\section{Minimal partitions and Courant-sharp}\label{smpCs}
The main object of this section is the proof of Theorem \ref{L=L}.
\subsection{Preliminaries}
Except possibly the question of eliminating the assumption
 that $\Omega$ is  simply connected, the proof in this section 
 uses essentially \cite{HH:2005a, HH:2005aprime} (or easy extensions of it). The minimal
 partitions which are involved
 in the proof are indeed regular. Although not very important here,
 this seems useful  to mention for possible extensions in higher
 dimensions
 where we do not have the fine results established in Section \ref{susanna2}.\\

\subsection{Definition of an exhausting family $N(u,\alpha)$.}
 Let $u$ be an  eigenfunction of $H(\Omega)$ with $k$         nodal domains
and  consider $N(u)\in \mathcal M(\Omega)$.
First we consider the finite sets of points
\begin{equation}\label{C*(N)}
C^*(N):=\mathcal Z_3 \cup (N(u)\cap \pa \Om).
\end{equation}
>From each of these points an arc emanates which ends either in the point itself (loop)
or ends in another point in $C^*(N)$. We call the collection of these arcs $\mathcal A_{*}$. 
Then we consider those components of $N(u)$ whose intersection with $C^*(N)$ is empty.
They have to be pairwise disjoint imbedded circles (without selfintersections) and we call the 
collection of these circles $\mathcal A_{**}$. Let us introduce
$$
\mathcal A = \mathcal A_{*}\cup \mathcal A_{**}\;.
$$
Note that each arc (or loop) $A\in \mathcal A$ is  rectifiable
 (because $N(u_m)$ is regular by Theorem \ref{reg_nodal})
and we can associate to $A$ naturally  a middle point $x_A$ in the
natural 
way ($x_A$ is chosen arbitrarily if $A\in \mathcal A_{**}$). 
We have a natural  arc length parametrization  starting from the point
$x_A$, but we prefer to parametrize $A$ as a parametrized curve 
 $[-1,+1]\ni t \mapsto L(A,t)$ such that $$L(A,0)=x_A\;,\;
 L(A,-1)=y_A^-\;, \;L(A,1)=y_A^+\;,$$
 where $y_{A^-}$ and $y_{A^+}$
 are the end points in the case of an arc, and where $y_{A^{-}}=y_{A^{+}}$ is the irregular point in the
 case
 of an irregular loop (i.e. in $\mathcal A_*$) 
and the opposite point in the case of a regular loop (i.e. in $\mathcal
A_{**}$).

For each $\alpha \in (0,1)$, we can consider the set 
\begin{equation}\label{Nt}
N(u,\alpha)=\{ N(u)\setminus  L(A,(-1+\alpha,1-\alpha ))\}
\end{equation}
and complete the definition by
 \begin{equation}
N(u,0)= \emptyset \mbox{ and }N(u,1)= N(u)\;.
\end{equation}
Note that by construction for every $0< \alpha $, $N(u,\alpha)$ contains all the critical 
points and  $N(u)\cap \pa\Om$; this will be important below.

\subsection{Proof of Theorem \ref{L=L}.}

{\bf We assume for contradiction that  for some $k$, $\mathfrak L_k=L_k$,
but that $\la_k<\la_m=L_k$ for some $m>k$.  }\\
Taking the smallest $m$ with this property, we can in addition assume
that 
\begin{equation}
\lambda_{m-1}<\lambda_m\;.
\end{equation}

Let $u_m$ a normalized eigenfunction such that $u_m$
 has $\mu(u_m)=k$ nodal domains $D_i$ ($i=1,2,\dots, k$).
 Then we associate with the exhausting family $N(u_m,\alpha)$ the
 decreasing 
 family
 of open sets~:
\begin{equation}
\Omega(\alpha) = \Om\setminus N(u_m,\alpha)\;.
\end{equation}
We want to consider the spectrum of~\\
\begin{equation}
H(\alpha):= H(\Om(\alpha))\,.
\end{equation}
 (We suppress the dependence on $\Om$ and $V$.)\\
Then,  $H(0)$ is our initial operator $H(\Omega)$ 
 and 
\begin{equation}
H(1)=\bigoplus_{i=1}^k H(D_i)
\end{equation}
Hence $H(1)$ has as lowest eigenvalue $\la_m$ with multiplicity $k$. 
By construction  $\sigma(H(0))=\sigma(H)$.
Furthermore $\la_1(H(1))$ has multiplicity $k$ and $$ \la_{k+1}(H(1))>\la_k(H(1))=\la_m(H(0))\;.$$
\begin{lemma}\label{monot}~\\
For any $\ell$, $\la_\ell(H(\alpha))$ is   monotonically increasing with $\alpha$.
\end{lemma}
{\bf Proof}~\\
 We just note that the  form domains $\mathcal Q(\alpha)$ of the 
quadratic forms $q(\alpha)$ associated to $H(\alpha),\;0\le \alpha\le 1$,
 satisfy 
 $\mathcal Q(\alpha)\subset\mathcal Q(\alpha')$ for $\alpha'\le \alpha$. 
    
\begin{lemma}\label{contal}~\\
For any $\ell$,  $\la_\ell(H(\alpha))$ depends continuously upon $\alpha$.
\end{lemma}
{\bf Proof}~\\
Although there is a lot of literature\footnote{We thank P.~Stollmann
  and M.~Dauge for useful discussions.} on the subject
(see \cite{Kato:1977,Simon:1977,Simon:1978,Weidmann:1980a, Weidmann:1980b, Weidmann:1984, Stollmann:1995}), it is difficult
 to give a reference corresponding to this crack situation.
This statement is proved  in Dauge-Helffer \cite{DH} (at least for the case
 of one crack). These authors treat the case when
 the condition on the crack is Neumann. The  Dirichlet case does not
 create new problems (the important point being the monotonicity which is evident
 in the case of Dirichlet). Note that we have strong convergence and
 that the left continuity and the right continuity should be treated
 separately.\\

We continue with   
\begin{lemma}\label{l==}~\\
 For 
each $\alpha\in [0,1]$, $\la_m\in\sigma(H(\alpha))$. 
\end{lemma}
\textbf{Proof} ~\\
By construction of $N(u_m,\alpha)$, the restriction of $u_m$
 to $\Omega(\alpha)$ is indeed an eigenfunction of 
$H(\alpha)$.
\finbox \\

We recall (see the two first lines of this subsection) that we are inside a proof by contradiction and consider
 first the following case.

\textbf{Case (a): $\mathbf{\la_m}$ is simple}\\
\begin{lemma}\label{casea}~\\
There is a minimal  $\alpha_1\in (0,1)$ such that 
\begin{equation}\label{salpha1}
\la_1(\alpha_1)<\la_2(\alpha_1)\le\dots\le\la_{m-1}(\alpha_1)=\la_m(\alpha_1)\,,
\end{equation}
with $\lambda_m(\alpha_1)=\lambda_m$.
\end{lemma}
\textbf{Proof}~\\
By assumption $m>k$ and $H(0)$ has $m$ eigenvalues smaller or equal to $\la_m$ whereas 
$H(1)$ has just $k$ eigenvalues smaller or equal to $\la_m$,   hence less. By continuity 
with respect to $\alpha$ at least one eigenvalue has to become larger than $\la_m$.
\finbox\\

Next we consider this eigenvalue of $H(\alpha_1)$. The restriction of
$u_m$
 to $\Omega(\alpha_1)$ gives a first eigenfunction and there exists  a
 real valued normalized  eigenfunction  $v$ of $H(\alpha_1)$ such that $v$
 is orthogonal to $u_m$ in $L^2(\Omega(\alpha_1))
$~:
\begin{equation}\label{orthv}
\langle u_m \;|\; v\rangle =0\;,
\end{equation}
where $\langle\,\cdot\;|\;\cdot\,\rangle$ denotes the $L^2$-scalar product.\\
We now play inside the two dimensional eigenspace spanned by $u_m$ and $v$.
\begin{lemma}\label{u,v}~\\
There exists $\beta_0>0$, such that 
$\forall \beta \in (-\beta_0,+\beta_0)$ 
the function $w_\beta=u_m+\beta v$ has exactly  $k$ nodal 
domains.\\
Furthermore, the family of the nodal domains of  $w_\beta$ gives, for $\beta \neq 0$,
a minimal bipartite  partition of $\Om$, which is distinct of the partition
 associated to $N(u_m)$.
\end{lemma}
\textbf{Proof}\\
We recall that this construction is done for $\alpha =\alpha_1\in(0,1)$.
For each $A \in \mathcal A$, let  
$$
I_A := L(A,(-1+\alpha_1,1-\alpha_1))\;,
$$
and let  $\mathcal V_A \subset \Omega $ be an open neighborhood  of 
$\overline{I_A}$, whose regular boundary
 crosses  $N(u_m)$ twice (transversally)
 and 
 such that each component of the open set  $\Omega (\alpha_1) \setminus \cup_{A\in
   \mathcal A} 
\overline{\mathcal V_A}$ (which is  contained in  $\Omega(1)$), is contained in a unique nodal domain of $u_m$.\\
We assume that we have colored these nodal domains by $+$ or $-$,
 and this permits to write,  for each $A\in \mathcal A$,  the decomposition
$$
\pa \mathcal V_A \cap \Omega(\alpha_1) = b_A^+ \cup b_A^-\;,
$$
where $b_A^\pm$
 is contained in a positive or negative nodal domain of $u_m$.\\

The first claim is now that there exists $\beta_0$ such that 
if we add $\beta v$, with $|\beta| \leq \beta_0$,  the number
 of nodal domains of $u_m +\beta v$ can only
 increase. Using  Hopf's boundary point lemma for $u=u_m$, 
\cite{GT:1983},  we have 
$|\nabla u_m(x)|>0$ for $x\in \left(N(u_m)\setminus
  N(u_m,\alpha_1)\right)\cap \Omega $ 
and using the property that $v$ vanishes at the boundary of $\Omega(\alpha)$, we obtain
the existence of $\beta_0$ such that $u_m+\beta v$ is strictly positive
 on each  $b_A^+$ and strictly negative on $b_A^-$.
It is then clear that associated to each positive $D_i$, there is
 at least one nodal domain of $u_m +\beta v$, with non trivial
 intersection with $D_i$ and contained in 
 $D_i \cup \left(\cup_A \mathcal V_A\right)$.
All these nodal domains are necessarily disjoint and this proves
 the first claim.\\

Let us now show that  
we cannot increase the number of  nodal domains. If it was
 the case, this would give an upper bound for $\mathfrak L_{k+1}$
 and using \eqref{frakL<}, we would obtain, using the strict monotonicity
 of the sequence $\mathfrak L_\ell$ (see \eqref{frakL<}) with respect to $\ell$, 
 $\lambda_m =\mathfrak L_k < \mathfrak L_{k+1}\leq \lambda_m$, 
 hence a contradiction.\\

So $w_\beta$ has also exactly $k$ nodal domains corresponding also to
a minimal  $k$-partition 
$\mathcal D'\in \mathfrak O_k$ of $\Omega$. But  $\mathcal D'\neq\mathcal D$ since both functions, $u_m$ and 
$v$ (and hence $w_\beta$ for $\beta\neq 0$) are linearly
independent. In addition, one can verify that $G(\mathcal D')$ is
bipartite.
 \finbox\\

Now we can complete the proof of Theorem \ref{L=L} for \textbf{case (a)}.\\
Indeed by Theorem \ref{partnod},
$w_\beta$ (more precisely the natural extension of $w_\beta$ to
$\Omega$) is an eigenfunction of $H=H(0)$ and 
therefore this would imply that $\la_m$ has multiplicity at least two, 
contradicting 
our assumption that $\la_m$ is simple.
\\
\begin{remark}~\\
Note that for these proofs we only need weak versions of our results
 because we work only with strong regular partitions (satisfying in
 addition the equal angle meeting property). So the techniques
 of \cite{HH:2005a} (as recalled in appendix \ref{sregcase}) are also relevant.

In the non simply connected case, if one wants to apply 
 \cite{HH:2005a}, one should also verify 
 a global compatibility condition for each homotopy class of $\Omega$. Because $w_\beta$ is 
an eigenfunction in $\Omega(\alpha_1)$, this is  a consequence
 of the property that any path in $\Omega$
 is homotopic to a path in $\Omega(\alpha_1)$. It is then easy
 to verify this  additional cycle-compatibility condition
 introduced in \cite{HH:2005a},  because $w_\beta$ is an
 eigenfunction of $H(\alpha_1)$.
\end{remark}

\noindent
\textbf{Case (b): $\mathbf{\la_m}$ has multiplicity greater than one.}\\
Assume that $\la_m$ has multiplicity $\ell >1$, so that the $m+\ell$
 first eigenvalues of $H=H(0)$ satisfy~:
\begin{equation}\label{sH0}
\la_1<\la_2\le\dots\le \la_k\le \dots\lambda_{m-1}<\la_m=\dots 
=\la_{m+\ell-1}<\la_{m+\ell}\;.
\end{equation}
Of course, all the previous constructions can be done but arriving
 at the last line of case (a), we loose the contradiction. The idea is that we have to choose our $w$ more carefully.\\

$\lambda_m$ being an eigenvalue of $H(\alpha)$ for any $\alpha$, we
can associate to $\lambda_m$ the 
eigenspace $\underline{U}(\alpha,\lambda_m)$ which is defined
 as the subspace in the spectral space of $H(\alpha)$ $U(\alpha,\lambda_m)$
 consisting of functions which are  restrictions to $\Omega(\alpha)$
 of eigenfunctions of $H(\Omega)$. Of course $\underline{U}(\alpha,\lambda_m)$
 contains $u_m$ but could be larger.\\
We then need to show that
\begin{lemma}\label{previouslemma}~\\
For $\alpha =\alpha_1$, $\underline{U}(\alpha_1,\lambda_m)$
 is strictly included in $ U(\alpha_1,\lambda_m)$.
\end{lemma}
{\bf Proof}~\\
For any $\alpha <\alpha_1$, we can choose a normalized 
$v_\alpha$ in $U(\alpha,\lambda_{m-1}(\alpha))$.
Then it is clear, observing that $\lambda_{m-1}(\alpha)<\lambda_m$,  that
\begin{itemize}
\item
$v_\alpha$ is orthogonal in $L^2(\Omega(\alpha))$ to $\underline{U}(\alpha_1,\lambda_m)$
 which is (more precisely, can be identified to) a subspace of
$U(\alpha,\lambda_{m})$ for any $\alpha < \alpha_1$.
\item
$v_\alpha$ is bounded independently of $\alpha$ in $H^1(\Omega(\alpha_1))$.
\end{itemize}
Then we can by compactness, find a sequence $w_{\alpha(n)}$
 such that $\alpha (n)$ tends to $\alpha_1$
 as $n \rightarrow + \infty$ and $v_{\alpha(n)}$
 converges weakly to some $v_{\alpha_1}$ in $W^1(\Omega(\alpha_1))$
 and strongly in $W^{s}(\Omega(\alpha_1))$ for $s<1$
 by compactness.\\
Now, it is clear that
\begin{itemize}
\item $v_{\alpha_1}$ is orthogonal to $\underline{U}(\alpha_1,\lambda_m)$.
\item $ ||v_{\alpha_1} ||=1$.
\item $ (-\Delta + V) v_{\alpha_1} =
 \lambda_m  v_{\alpha_1}\mbox{ in }\Omega(\alpha_1)\;.$
\end{itemize}
With a small additional work, one can show 
 that $v_{\alpha_1}\in W_0^1(\Omega(\alpha_1))$.
So  $v_{\alpha_1}$ is effectively in the
 form domain for  the Dirichlet problem  in $\Omega(\alpha_1)$ and  in
 $U(\alpha_1,\lambda_m)\cap \underline{U}(\alpha_1,\lambda_m)^{\perp}$.
\paragraph{End of the proof of case (b).}~\\

The argument is then as in case (a), but, using Lemma \ref{previouslemma},
  we can choose a non trivial $v$
 in $U(\alpha_1,\lambda_m)\setminus\underline{U}(\alpha_1,\lambda_m)$.
But on one  hand $w_\beta $ cannot belong to 
$\underline{U}(\alpha_1,\lambda_m)$ (because $v$ does'nt).
On the other hand, we have obtained some  $\beta\neq 0$ such that $w_\beta$ extends as an eigenfunction  of $H (\Omega)$
 hence by definition  in $\underline{U}(\alpha_1,\lambda_m)$
 and this gives the contradiction.

Theorem \ref{L=L} has an immediate consequence.
\begin{corollary}\label{T+Pl}~\\
Let $\Om$ satisfy Assumption \ref{Asswhole1} and assume that $V\in L^\infty(\Omega)$. Then 
\begin{equation}\label{Lk=Lk}
\big|\{k\;|\;\mathfrak L_k=L_k\}|\;<\;\infty.
\end{equation}
\end{corollary}
\textbf{Proof}~\\
This is an immediate consequence of Theorem \ref{L=L} and of Pleijel's 
Theorem~\ref{Pleijel}. \finbox
 
\begin{remark} ~\\
It is now easier to analyze the situation for the disk and  for rectangles
 (at least in the irrational case), since we have  just to check for which eigenvalues
 one can find associated Courant-sharp eigenfunctions.  This will be
 done in  sections \ref{sdisk} and \ref{srectangle}.
\end{remark}
\newpage
\section{Further properties of subpartitions}\label{ssubp}
All the statements of this section illustrate the rigidity of 
 the structure of the subpartitions. This can be very efficient 
 for disproving that a partition is minimal. In particular we will
 prove Theorem \ref{thsubpartition}.\\

The following proposition is useful :
\begin{proposition}\label{subpartition}~\\
Under Assumptions \ref{Asswhole1} and \ref{Asswhole2},  
let $\mathcal D = (D_i)_{i=(1,\dots,k)}$ be  a $k$-minimal partition
 for $\mathfrak L_k(\Omega)$. Then, for any subset $I\in
 \{1,\dots,k\}$, the associated subpartition  $\mathcal D^I = (D_i)_{i\in I}$ satisfies 
\begin{equation}
\mathfrak L_k= \Lambda (\mathcal D^I) = \mathfrak L _{|I|} ( \Omega^I)\;,
\end{equation}
 where 
$$
\Omega^I:=\Inte (\overline{\cup_{i\in I}  D_i})\;.
$$
\end{proposition}

{\bf Proof}~\\
We prove this proposition by contradiction. If it was not the case, we
would
 construct (starting of a minimal $|I|$-partition of $\Omega^I$) a
new  minimal partition $\widetilde D$  of $\Omega$, for which the
 $\lambda(\widetilde D_i)$'s are not equal in contradiction
 with what we proved in Section \ref{susanna2} (see also (d)
 in Remark \ref{osservazioni} together with the fact that
 $k_0=\emptyset$
 in our case).\\

 As a consequence of a
more general theorem in \cite{AHH:2004},  we have (with very weak
assumptions on $\Omega$)
 the analogous of Proposition \ref{subpartition}
\begin{proposition}\label{Csharp}~\\
Suppose $u$ is \textbf{Courant-sharp}. Denote the associated  nodal domains
by $\{D_i\}_1^k$. Let $L$ be a subset of $\{1,2,\dots,k\}$ with 
$\# L=\ell < k$ and let $\Omega^L=\Inte(\overline{\cup_{i\in L}D_i})\setminus \pa \Om$. Then 
\begin{equation}\label{lal}
\la_\ell(\Omega^L)=\la_k
\end{equation}
where $\la_j(\Omega^L)$ are the eigenvalues of $H(\Omega^L)$.

Moreover, if $\Om^L$ is connected, $u\big|_{\Om^L}$ is \textbf{Courant-sharp}
and $\la_\ell(\Om^L)$ is simple.
\end{proposition}

\begin{proposition}~\\
Under the assumptions of the introduction on $\Omega$ and $V$, let $\mathcal D=(D_1,\dots,
D_k)$ be a minimizing partition relative to $\mathfrak L_k$. Let
$\mathcal D^\prime\subset \mathcal D$ be any subpartition into $1\leq
k^\prime\leq k$ elements which is bipartite relatively to $$ \O^\prime:=\Inte\left(\bigcup_{D_i\in\mathcal D^\prime}\overline{D_i}\right)\;.$$
Then 
\begin{itemize}
\item[(a)] $\mathfrak L_k=\lambda_{k^\prime}(\O^\prime)$;
\item[(b)] If $k^\prime<k$,  and $\O'$ is connected 
  then $\lambda_{k^\prime}(\O^\prime)$ is simple.
\end{itemize}
\end{proposition}

{\bf Proof}~\\
The point (a) is proven like at the end of Section \ref{susanna2}. The
proof
 of (b) is an immediate consequence of Proposition \ref{Csharp}, if we 
add the assumption that there exists a  bipartite subpartition $\mathcal D''$
 such that  $ \mathcal D^\prime\subset \mathcal D''\subset  \mathcal
 D$ with  $k'<k''\leq k$.\\
The general proof is a little more tricky.
Given the minimal partition $\mathcal D$ of $\Omega$, and the
subpartition  $\mathcal D'$ (with associated bipartite graph),
 there is a
subfamily $\mathcal A_0$ of the set $\mathcal A$ 
 of arcs (we keep the ends of the arcs living
 in $\Omega$) belonging to $N(\mathcal D)$ (see the discussion
 in Section \ref{smpCs}) such that
$\widetilde \Omega:=\Omega \setminus \mathcal A_0$ has for the same
 partition  a bipartite graph relatively to $\widetilde \Omega$
 and  such that
$\mathcal A_0$ does not contain
  arcs belonging to the intersection of the boundaries
 of two neighbors  of $\mathcal D'$ in $\Omega'$.
 But we are in a Courant-sharp situation, so 
\begin{equation}\label{share}\mathfrak L_k(\Omega)
 =\mathfrak L_k(\widetilde \Omega) =\lambda_k(\widetilde \Omega)\;,
\end{equation} 
  and we can apply
 the Courant-sharp Theorem relative to $\widetilde \Omega$ and the
 partition $\mathcal D$ (which is now the family
 of nodal domains of an eigenfunction on $\widetilde \Omega$).\\
\begin{remark}~\\
In the preceding proof, $\widetilde \Omega$ is not unique, but
 it is interesting to emphasize that what we have shown is that, 
 once a minimal $k$-partition of $\Omega$ is given, then 
 all the possible $\widetilde \Omega$'s should share the property
 \eqref{share}.
\end{remark}

In order to complete the proof of Theorem \ref{thsubpartition}, it remains to
establish our uniqueness result at the level of the subpartitions
 of a minimal partition.

\begin{proposition}[Uniqueness]~\\
Let $\mathcal D$ be a minimal $k$-partition relative to $\mathfrak
L_k(\Om)$. Let $\mathcal D^\prime\subset \mathcal D$ be any
subpartition of $\mathcal D$ into $1\leq k^\prime< k$ elements and let
$$\O^\prime= \Inte\left(\bigcup_{D_i\in\mathcal
  D^\prime}\overline{D_i}\right)\;,$$
be connected. 
Then  $\mathfrak L_{k^\prime}(\O^\prime)$  is uniquely achieved.
\end{proposition}
We know already from Proposition \ref{subpartition} that $\mathfrak L_{k^\prime}(\O^\prime) = \mathfrak
L_k(\Om)$.
The proof of uniqueness is by contradiction. Let $I$ and $J$ two subsets of
$\{1,\dots,k\}$
 such that $I\cap J =\emptyset$ and $I\cup J = \{1,\dots,k\}$.\\
If we have indeed two minimal subpartitions of $\Omega_I$ for some $I$ of 
 cardinality strictly less than $k$. We can complete these two
 subpartitions by the open sets $D_j$ ($j\in J$). According to the
 proof of our results. Now take a pair $(D_i,D_j)$, with $i\sim j$
 ($i\in I$
 and $j\in J$). This should exist if $\Omega$ is connected and
 $\Omega_I$ is connected. There is necessarily another $\tilde
 D_{\tilde i}$ of the other partition which meets $D_i$
 and is a neighborhood of $D_j$. But the eigenfunction
$u_{ij}$ of $D_{ij}$ and $\tilde u_{ij}$ of $\tilde D_{ij}$
 should be proportional on $D_j$ to $u_j$. It is then clear that by
 unique
 continuation  $D_i$ should coincide with $\tilde D_i$
 and $u_i$ should be proportional to $\tilde u_i$.
Possibly iterating the argument, we arrive to a contradiction that the
two partitions are different.
\newpage 

\section{Example 1: The case of the disk}\label{sdisk}
In this section, we analyze in great detail the case of the disk. Although the spectrum is explicitly computable, we are mainly interested in the ordering of the eigenvalues corresponding to different angular momenta. In particular this will
 give   a first example where the inequalities in (\ref{ineqbetw}) can be strict.

Consider the Dirichlet realization $H_0$ in  the unit disk $B_1\subset\mathbb R^2$. 
We have in   polar coordinates~:
$$
-\Delta=- \frac{\pa^2}{\pa r^2} -\frac 1r \frac {\pa}{\pa r} 
- \frac{1}{r^2} \frac{\pa^2}{\pa \theta^2}\;,
$$
and the Dirichlet boundary conditions require that any eigenfunction $u$ 
satisfies $u(r,\theta)=0$ for $r=1$.
We analyze for any $\ell\in \mathbb N$ the eigenvalues
 $\lambda_{\ell, j}$ of 
$$
(- \frac{d^2}{d r^2} - \frac 1r \frac {d}{dr} + \frac{\ell^2}{r^2})f_{\ell,j}=
\la_{\ell,j}f_{\ell,j}\;,\; \mbox{ in } (0,1)\;.
$$
We observe that the operator is self adjoint for the scalar product in
 $L^2((0,1),r\,dr)$.\\
The corresponding eigenfunctions of the eigenvalue  problem take the
 form 
\begin{equation}\label{diskf}
u (r,\theta)=f_{\ell,j} (r) \left( a\cos\,\ell \theta\; + \;
b \sin \,\ell \theta\,\right),\:\mbox{ with } a^2+b^2>0\;,
\end{equation} 
where the $f_{\ell,j}(r)$ are suitable Bessel functions satisfying 
for $\ell=0$,  $f_{0,j}'(0)=0$ and $f_{0,j}(1)=0$ and for $\ell>0$, $f_{\ell,j}(0)=
f_{\ell,j}(1)=0$.
For  the corresponding $\lambda_{\ell, j}$'s,  we find by looking 
up for instance 
\cite{Bandle:1980} (or in appendix \ref{sBessel})  the following ordering.
\begin{equation}\label{ladisk}
\begin{array}{l}
\la_1=\la_{0,1}<\la_2=\la_3=\la_{1,1}<\la_4=\la_5=\la_{2,1}<\la_6=\la_{0,2}<
 \lambda_7 = \lambda_8 = \lambda_{3,1} <\dots\\
\qquad \dots < \lambda_9=\lambda_{10} =\lambda_{1,2} <
 \lambda_{11}=\lambda_{12} =\lambda_{4,1} < \dots \\
\qquad \dots < \lambda_{13}=\lambda_{14} = \lambda_{2,2} <
 \lambda_{15}=\lambda_{0,3}<\dots.
\end{array}
\end{equation}
 We recall that the zeros of the Bessel functions are related  to the eigenvalues by the relation
\begin{equation}\label{zeroeig}
\lambda_{\ell,k} = \,( j_{\ell,k})^2\;.
\end{equation}

We hence have from \eqref{diskf}
\begin{equation}\label{mdisk}
\begin{array}{l}
\mu(u_1)=1,\:\\
 \mu(u)=2, \mbox{ for any eigenfunction u associated to } \lambda_2=\lambda_3,\\
\mu(u)=4,  \mbox{ for any eigenfunction u associated to } \lambda_4=\lambda_5,\\
 \mu(u_6)=2,\\
\mu (u) =6,  \mbox{ for any eigenfunction u associated to }
 \lambda_7=\lambda_8, \\
\mu (u) = 4,  \mbox{ for any eigenfunction u associated to }
 \lambda_9=\lambda_{10},\\
 \mu (u) = 8,  \mbox{ for any eigenfunction u associated to }
 \lambda_{11}=\lambda_{12},\\
\mu (u) = 8,  \mbox{ for any eigenfunction u associated to }
 \lambda_{13}=\lambda_{14},\\
\mu(u_{15}) =3
\;.
\end{array}
\end{equation}
Hence 
\begin{equation}\label{L1lambda}
L_1=\lambda_1\;,\; L_2=\lambda_2\;,\;L_3=\lambda_{15}\;, \; L_4=\lambda_4
\end{equation}
and this implies 
\begin{equation}\label{L3lambda}
L_3 > \lambda_6 > L_4 > \lambda_3\;.
\end{equation}
 In addition, let us show that 
\begin{equation}\label{L3L3}
 L_3>\mathfrak L_3.
\end{equation}

This can be seen as follows. We can split $B_1$ in three  
sectors with opening angle
$2\pi/3$. Call such a sector $S_{1/3}$ then the corresponding 
eigenvalue $\lambda(S_{\frac 13})$
 satisfies~:
$$
\mathfrak L_3 \leq \lambda(S_{\frac 13})
$$
We then observe that by monotonicity that~:
$$
\lambda(S_{\frac 13}) < \lambda(S_{\frac 14})$$
and we can recognize that  $\lambda(S_{\frac 14})=\lambda_4$ observing that $\lambda_4 = \lambda_{\ell,j}$
 with $\ell =2$ and $j=1$ (see (\ref{diskf})). The proof of (\ref{L3L3}) is then
 a consequence of (\ref{L3lambda}).
\begin{remark}~\\
An open problem is to prove or disprove the equality $\mathfrak L_3 = \lambda(S_{\frac 13})$.\\
A possibility for trying to find other configurations could be
 to look at the second eigenvalue of  a half disk problem 
$\{(x_1^2+x_2^2) < 1\}\cap \{x_2 >0\}$, where we 
 take Dirichlet on the circular part and a part of the straight basis
 and Neumann on the other part. More precisely, we take Neumann
 on $x_2=0\;,\; |x_1|<t$, where $t$ is a free parameter. The nodal domains  of the second eigenfunction (completed by symmetry) could give an alternative  candidate for such a minimal partition.\end{remark}

In the case of the disk, we have
\begin{proposition}~\\
Except the cases $k=1$, $2$ and $4$, minimal partitions
 never correspond to nodal domains.
\end{proposition}
{\bf Proof}~\\

According to Theorem \ref{L=L}, it is enough to investigate when the $k$-th
 eigenvalue corresponds to  Courant-sharp eigenfunctions.

One can in addition use the twisting trick (see \cite{HH:2005b}) for
 eliminating all the eigenvalues $\lambda_{\ell,m}$, for which
 $m\geq 2$ and $\ell >0$. This trick goes roughly as follows. When
 $\ell >0$, we can divide the disk as the union of a smaller disk
 and of its complementary, each of these sets being the union of at
 least
 two nodal domains. Then by small rotation of the small disk, we obtain
 a new partition which has the same energy. If the initial one was
 minimal, the new one should be also minimal, but it is easy to show
 that
 the new one has not the ``equal angle meeting'' property 
 of a regular partition. This gives
 the contradiction.\\

So we have finally to analyze the eigenvalues $\lambda_{0,k}$
 and the family $\lambda_{\ell,1}$.

For the first family, we observe that $\lambda_{0,k}$
 can neither be the $k$-th eigenvalue as soon as $k\geq 2$.

For the second family, which occurs only for $k=2\ell$ even,
 inspection of the tables leads to the condition $k\leq 4$, we observe indeed
 that $\lambda_{0,2} < \lambda_{3,1}$.

We also observe that when $k$ is odd, we obtain that necessarily $L_k=\lambda_{0,k}$.

\begin{remark}~\\
It could be interesting to determine when $\mathfrak L_k$
 is an eigenvalue of the Laplacian on the double covering.
Again, one can show (see \cite{HH:2005b}) that this cannot
 be the case for $k$ large.
\end{remark}
\newpage 
\section{Example 2: the case of the rectangle}\label{srectangle}
Note that for the case of a rectangle, 
the spectrum and the properties of the eigenfunctions are analyzed
 as toy models  in
 \cite{Pleijel:1956}, Section 4. This was also used for testing general conjectures in \cite{AHH:2004}.\\
For a rectangle of sizes $a$ and $b$, the spectrum is given by
$\pi^2 (m^2/a^2 + n^2 /b^2)$ ($(m,n)\in (\mathbb N^*)^2$).\\
The first remark is that all the eigenvalues are simple
 if $\frac{a^2}{b^2}$ is irrational. 
Except for specific  remarks for the square, we now assume\\
$$
 (a/b)^2 \mbox{ is  irrational}.
$$
 So we can associate to each eigenvalue
$\lambda_{m,n}$, an (essentially) unique eigenfunction $u_{m,n}$
 such that $\mu(u_{m,n}) = nm$.
Given $k\in \mathbb N^*$, the lowest eigenvalue corresponding to $k$ nodal domains is given by
$$
L_k = \pi^2 \inf_{mn=k} (m^2/a^2 + n^2 /b^2)\;.
$$

The behavior of $L_k$ can depend dramatically of the arithmetical properties
 of $k$ but what is important for us is that in any case we have
\begin{equation}\label{last1}
L_k \geq 2 \pi^2 k/ (ab)\;.
\end{equation}
This immediately implies
\begin{equation}\label{last2}
\liminf_{k\ar +\infty} \frac{L_k}{k} \geq \frac{2 \pi^2 }{ab}\;.
\end{equation}
We note that the right hand side can also be written in the form
$\lambda([0,1]^2)/|\Omega|$.\\

\begin{remark}~\\
Note that these estimates are much better than the estimates obtained by
 Faber-Krahn inequality (see in Section \ref{spa}).
\end{remark}

As we have seem in  \eqref{comphexa} and using\footnote{ We thank
 V. Bonnaillie-No$\ddot{e}$l and G. Vial for giving us a precise numerical approximation of $\lambda(Hx1)$. According to A.~El Soufi, it seems unknown
 that the hexagon gives the minimal eigenvalue between all the
 polygons (of same area) permitting to realize a perfect partition  of the plane. We just
 compare here the square and the hexagon.} the comparison between the lowest eigenvalue of the hexagon and  of the square,
\begin{equation}
\lambda(Hx1) \sim 18.59013 < \lambda ([0,1]^2)= 2 \pi^2\sim 19.7392\;,
\end{equation}
 this will imply

\begin{equation}
\liminf_{k\ar +\infty} \frac{L_k}{k} > \limsup_{k\ar +\infty} \frac{\mathfrak L_k}{k}\;.
\end{equation}
This implies that $L_k>\mathfrak L_k$ for $k$  large.

\begin{remark}~\\
In the case when $(\frac ab)^2$ is rational we could have problems in
 the case of multiplicities. We have then to 
 control the nodal sets of the eigenfunctions corresponding to the
 degenerate eigenvalues which are  $\leq \mathfrak L_k$.
\end{remark}

We now describe all the possible situations.\\
\begin{lemma}~\\
In the irrational case, $\lambda_{m,n}$ cannot lead to 
 a Courant-sharp situation
 if\break   $\inf (m,n)\geq 3$.
\end{lemma}

{\bf Proof}~\\
Applying Proposition \ref{Csharp},  it is sufficient to analyze the case when $m=n=3$.
It is then enough to show that $\lambda_{3,3}$ cannot be the ninth eigenvalue.\\
Because the eigenvalues corresponding to $\max(m,n)\leq 3$
 are obviously below,  let us assume by contradiction
 that $\lambda_{3,3} \leq \lambda_{1,4}$ and that $\lambda_{3,3} \leq 
\lambda_{4,1}$.

This reads
$$
\frac {9}{a^2} + \frac{9}{b^2} \leq \frac {1}{a^2} + \frac{16}{b^2}\;,
$$
and
$$
\frac {9}{a^2} + \frac{9}{b^2} \leq \frac {1}{a^2} + \frac{16}{b^2}\;.
$$

So, we obtain
$$
\frac 87 \leq \frac{a^2}{b^2} \leq \frac 78\;,
$$
Hence a contradiction.

The next step is given in 
\begin{lemma}~\\
In the irrational case, $\lambda_{m,n}$ cannot lead to a  Courant-sharp situation
 if $m=2$ and $n\geq 4$ or if $m\geq 4$ and $n=2$.
\end{lemma}
Again it is enough to look at the case $(m=2,n=4)$, and to show that it cannot be the eight-th eigenvalue. Using the same idea as in the previous
 lemma, we assume by contradiction that
 $\lambda_{2,4} \leq \lambda_{1,5}$ and that $\lambda_{2,4} \leq 
\lambda_{3,1}$.\\
This reads
$$
\frac {4}{a^2} + \frac{16}{b^2} \leq \frac {1}{a^2} + \frac{25}{b^2}\;,
$$
and
$$
\frac {4}{a^2} + \frac{16}{b^2} \leq \frac {9}{a^2} + \frac{1}{b^2}\;.
$$
So, we obtain
$$
\frac 13 \leq \frac{a^2}{b^2} \leq \frac {1}{3}\;.
$$
But this gives $\frac{a^2}{b^2}=\frac 13$, which is excluded by
 the assumption that $\frac{a^2}{b^2}$ is irrational.\\
Let us now analyze the remaining cases.\\

\noindent {\bf In the case $m=2$, $n=3$}, an eigenfunction corresponding to
$\lambda_{m,n}$ is 
Courant-sharp if 
$$
\frac {4}{a^2} + \frac{9}{b^2} \leq \frac {9}{a^2} + \frac{1}{b^2}\;,
$$
and
$$
\frac {4}{a^2} + \frac{9}{b^2} \leq \frac {1}{a^2} + \frac{16}{b^2}\;.
$$
So, we obtain
$$
\frac 8 5 \leq \frac{a^2}{b^2} \leq \frac 53 \;.
$$

\noindent {\bf The case $m=3$, $n=2$,} is obtained by exchanging the role of $a$ and $b$.
So, we obtain
$$
\frac 8 5 \leq \frac{b^2}{a^2} \leq \frac 53 \;.
$$

\noindent {\bf In the case $m=2$, $n=2$,} we obtain similarly
$$
\frac 35 \leq \frac{a^2}{b^2} \leq \frac 53 \;.
$$

\noindent {\bf For the case $m=1$, $n=k$,} we obtain
$$
\frac {1}{a^2} + \frac{k^2}{b^2} \leq \frac {4}{a^2} + \frac{1}{b^2}\;,
$$
So, we obtain simply
$$
\frac{(k^2-1 )}3 < \frac{a^2}{b^2} \;.
$$

\noindent Finally, {\bf the case $m=k$, $n=1$} leads to the condition 
$$
\frac{(k^2-1 )}3 < \frac{b^2}{a^2} \;.
$$

\paragraph{A candidate for the $3$-minimal partition on the
  square.}~\\

In the case of the square, an argument similar to the case of the disk shows
 that $\mathfrak L_3$ (which should be smaller than $\mathfrak L_4$)
 is strictly less than $L_3$. We observe indeed  that $\lambda_4$
 is Courant-sharp, so $\mathfrak L_4 =\lambda_4$, and there is no
 eigenfunction corresponding to $\lambda_2=\lambda_3$
 with three nodal domains (by Courant's Theorem).

 \begin{figure}[h!]\label{picturesquare}
\begin{center}
\includegraphics[width=12cm]{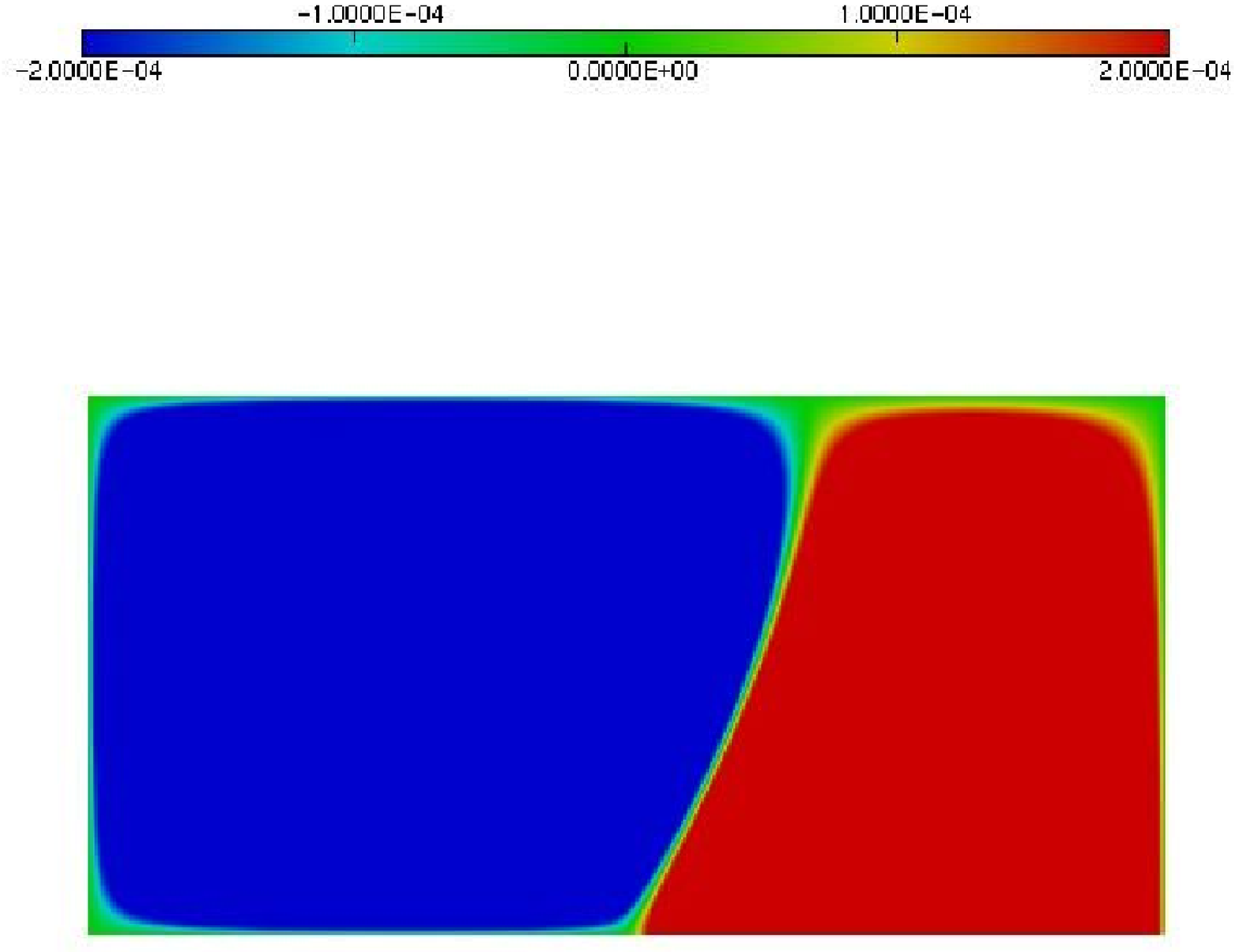}
 \end{center}
\end{figure}
Restricting to the half-rectangle and assuming that
 there is a minimal partition which is symmetric with one of the symmetry
 axes of the square perpendicular to two opposite sides, one is reduced to analyze a family of
 Dirichlet-Neumann
 problems. Numerical computations\footnote{ See 
 http://www.bretagne.ens-cachan.fr/math/people/virginie.bonnaillie/BH/MinimalPartitions/}  performed by V.~Bonnaillie-Noel (in
 January 2006) and G.~Vial 
 lead to a natural candidate (See Picture above) for a symmetric minimal partition.

The complete structure is reobtained from the half square 
by symmetry with respect to the horizontal axis. We observe
numerically
 that the three lines of $N(\mathcal D)$ meet at the center of the
 square.
 As expected by the theory they meet with equal angle $\frac{2\pi}{3}$
 and start from the boundary orthogonally.\\

\newpage 
\appendix
\section{Zeros of Bessel functions}\label{sBessel}
Let $j_{\ell,k}$ the $k$-th zero of the Bessel function corresponding to
 the integer $\ell\in \mathbb N$. The reference
 is the book by G.N. Watson \cite{Wa}. The most important statement for
 us is that $j_{\ell,k} =j_{\ell',k'}$
 imply if $\ell$ and $\ell'$ are positive integers,
 that $\ell =\ell'$ and $k=k'$. 
We refer to the subsection 15.28 (p. 484-485). Note that the proof of this result
 is based on deep results by Siegel about algebraic numbers.\\

Here is a list of approximate values
 after the celebrated handbook of \cite{AS}, p. 409, we keep only the values which are less than approximately $13$.

\begin{equation}
\begin{array}{rcccccccccccc}
&\ell=&0 &1 &2 & 3 & 4 & 5 & 6&7&8\\
k=
 1 && 2.40& 3.83 & 5.14 &6.38&7.59& 8.77&9.93&11.08&12.22&\\
2& & 5. 52 & 7.02 & 8. 42 & 9.76 & 11. 06 & 12.34&.&&\\
3& &8,65&10.17&11.62 & 13.02&..&&&&\\
4& & 11.79&13.32&.&&&&&&
  \end{array}
\end{equation}
This leads to the following ordering of the zeros :
\begin{equation}
\begin{array}{l}
j_{0,1}< j_{1,1} <
 j_{2,1} < j_{0,2}< j_{3,1}<j_{1,2}<j_{4,1}<j_{2,2}<j_{0,3}<\dots\\
\qquad \dots  < j_{5,1} < j_{3,2} < j_{6,1} < j_{1,3} < j_{7,1} < j_{2,3}
 < j_{0,4} < j_{8,1}\;.
\end{array}
\end{equation}

Note that, using Sturm-Liouville theory, (see for
 example the proof in \cite{Wa}),  the
following inequalities are always true~:
\begin{equation}\label{inequtz}
j_{\ell,k} < j_{\ell +1,k} < j_{\ell, k+1}\;,\; \forall \ell \in
\mathbb R^+\,,
\forall k \in \mathbb N^*\;.
\end{equation}
As a corollary, we obtain
\begin{equation}
j_{\ell,k} \leq j_{0, k + \ell}\;,
\end{equation}
with strict inequality for $\ell >0$.\\
It is also useful to have the half integer results (see in \cite{AS}, p. 467)
\begin{equation}
\begin{array}{rcccccccc}
&\ell=&\frac 12 & \frac 32 & \frac 52 & \frac 72 & \frac 92 & \frac {11}{2}&\frac{13}{2} \\
~&\\
k=
 1 && \pi& 4.49 & 5.76&6.99&8.18& 9.36& 10.51 \\
2& & 2 \pi & 7.73 & 9.09 & 10.42 & 11.7&\\
3& & 3 \pi&10.9&12.52 & 13.69&..\\
4& & 4 \pi & 14.06&
\end{array}
\end{equation}

\section{Alternative approach in the regular case}\label{sregcase}
Although not needed in this article, we recall some of the results of
 Statements of \cite{HH:2005a, HH:2005aprime}. 
The  main theorem is the following
  
\begin{theorem}\label{specconv}~\\
Suppose that  $\Om$ is regular and simply connected and that,
 for some  regular closed set $N$  satisfying in addition the ``equal angle
 meeting condition''. Suppose that, for some 
   $\la\in \mathbb R$, the associated family 
$\mathcal D(N)=\{D_1,\dots,D_\mu\}$ is admissible and 
satisfies a Pair Compatibility Condition, which means that $\lambda$
 is an eigenvalue of $H(D_{i,j})$ for which $D_i$ and $D_j$ are the
 two 
 nodal domains of some corresponding eigenfunction. Then there 
is an  eigenfunction $u$  of  $H(\Om)$ with corresponding eigenvalue  $\la
 $ such that  the family  of nodal domains of $u$ is $\mathcal D(N)$.
\end{theorem}
\begin{remarks}~
\begin{itemize}
\item 
If $\Om$ is not simply connected then the result does not hold in
general. One should add
 a nonholonomy condition (see \cite{HH:2005aprime}). In the case of
 minimal partitions, we have seen that this condition is reduced
 to the bipartite condition.
\item 
Note that the Pair Compatibility Condition is weaker than to assume
that
 $\lambda$ is the second eigenvalue of $H(D_{i,j})$ for each pair of  neighbors
 $(D_i,D_j)$.
\end{itemize}
\end{remarks}

As an application of this theorem, the authors  obtain~:
\begin{corollary}\label{Appli}~\\
Let $\Omega$ be  simply connected, $k\in \mathbb N$ ($k\geq 2$) and let 
$\mathcal D^{min}=(D_i)_{i=1,\dots,k}$ be a
 minimal  admissible strong regular\footnote{The notion of regularity
   was
 actually stronger there, but the Pair Compatibility Condition
 gives actually some regularity assumption on the boundaries.} 
partition. Then 
 there is an eigenfunction $u$
 of $H(\Omega)$ associated with
$$
\lambda = \max_i (\lambda(D_i))\;,
$$
such that $\mathcal D^{min}$ is the family
 of the $k$ nodal domains of $u$.
\end{corollary}
{\bf Proof}\\
Let us apply  Theorem \ref{specconv}. 
 We take as $\lambda  = \max_i (\lambda(D_i))$.\\
  The first point is 
 that all the $\lambda(D_i)$ should be equal. If not, one could by deformation
 of the $D_i$'s in a neighborhood of  regular points of their boundary find a new partition $\widetilde{\mathcal D}$,
 which would decrease $\max_i (\lambda(D_i))$.\\
The second point is to observe that considering two neighbors $D_i$ and $D_j$,
 then $\lambda$ should be the second eigenvalue of $H(D_{i,j})$. If it was not the case for some pair $(i,j)$, the two nodal domains of the second eigenfunction of $H(D_{i,j})$
 will give two new open sets $D'_{i}$ and $D'_j$
 with $\lambda(D'_i)=\lambda(D'_j)$, in contradiction with the 
 assumption of minimality and the first point  of the proof.\\
Hence the Pair Compatibility Condition is satisfied.\\

\begin{remark}~\\
As mentioned in the introduction, the case where $k=2$ corresponds to a rather well known characterization of the second eigenvalue of $H(\Omega)$. The admissibility condition is
 of course automatically satisfied in this case.
\end{remark}

\newpage

\footnotesize
\bibliographystyle{plain}

\scshape 
B. Helffer: D\'epartement de Math\'ematiques, Bat. 425,
Universit\'e Paris-Sud, 91 405 Orsay Cedex, France.

email: Bernard.Helffer@math.u-psud.fr

\scshape
T. Hoffmann-Ostenhof: Institut f\"ur Theoretische Chemie, Universit\"at
Wien, W\"ahringer Strasse 17, A-1090 Wien, Austria and International Erwin
Schr\"odinger Institute for Mathematical Physics, Boltzmanngasse 9, A-1090
Wien, Austria.

email:thoffman@esi.ac.at

\scshape
S. Terracini: 
Universit\`a di Milano Bicocca, Via Cozzi, 53 20125 Milano (Italy).

email: susanna.terracini@unimib.it

\end{document}

%% file: susanna06m.tex












\def\pt#1{{\bf Proof of Theorem \ref{#1}.}}
\def\eq#1{(\ref{#1})}

\def\neweq{\begin{equation}}
\def\endeq{\end{equation}}

\def\O{\Omega}
\def\bbbr{\mathbb{R}}

\section{Optimal partitions in $N$--dimensions}\label{susanna1}
Let $\Omega\subset\bbbr^N$ be a connected, open bounded domain. Note
that 
at this stage we don't need any regularity of the boundary. 
However we will need it later (in Section \ref{susanna2}), in the $2$-dimension
case 
 in order to describe the 
local structure of nodal lines at their intersection with the boundary.

\begin{definition}\label{eigenvalue}~\\
For any measurable $D\subset \Omega$ and for $V\in L^\infty(\Omega)$, let $\lambda_1(D)$ denotes the
first eigenvalue of the Dirichlet realization of the  Schr\"odinger operator in the following generalized sense. 
We define
$$
\lambda_1(D)=+\infty\;,
$$
if
$
\left\{ u\in W^1_0(\Omega)\,, u\equiv 0\; \text{a.e. on}\;
 \Omega\setminus D \right\}=\left\{0\right\}
$, 
and
$$
\lambda_1(D)=\min\left\{\frac{\int_\Omega \left(|\nabla
u(x)|^2+V(x)u(x)^2\right)dx}{\int_\Omega|u(x)|^2dx}\;:\;
 u\in W^1_0(\Omega)\setminus\{0\}\,, u\equiv 0\; \text{a.e. on}\;
 \Omega\setminus D \right\}\;,
 $$
otherwise. We call groundstate any function $\phi$ achieving the above minimum. 
\end{definition}

We shall always assume that
$$
\lambda_1(\O)>0\;.
$$

\begin{remark}\label{rempotential}~\\
The presence of an  $L^\infty$ potential $V$ does not create
 particular problems. We prefer to simplify the notations
 to  explain all the proofs with the additional assumption that $V$ is
 identically $0$. In this case the positivity on $\lambda_1(\Omega)$
 is effectively satisfied. In the general case, we can always assume
 this property
 by adding a constant to $V$.
\end{remark}
We observe  that the minimization problem always possesses a (possibly not unique) non negative solution $\phi\geq 0$. We shall always make this choice. Next we consider the following class of minimal partition problems~:
\neweq\label{model}
\mathfrak L_{k,p}:=
\inf_{{\mathcal B}_k}\left(\frac{1}{k}
\sum_{i=1}^k \big(\lambda_1(D_i)\big)^p\right)^{1/p}\,,
\endeq

\neweq\label{modelbis}
\mathfrak L_k:=
\inf_{{\mathcal B}_k}\max_{i=1,\dots,k} \big(\lambda_1(D_i)\big)
\endeq
where the minimization is taken over the class
 of partitions in $k$ ``disjoint'' 
measurable subsets of $\O$
$${\mathcal B}_k:=\left\{\mathcal{D}=(D_1,\dots,D_k):
\;\bigcup_{i=1}^k{D}_i
\subset\O,\;
 |D_i\cap D_j|=0\;{\rm if}\;i\neq j\right\}\;.$$

\begin{remark}\label{relaxation}~\\
 The values $\mathfrak  L_{k}$  considered in this section can be viewed as a relaxation
of those defined in the introduction. We have indeed replaced
 ``open'' by ``measurable''. We keep the same notation, for we shall prove as a part of our regularity theory that, in all the
interesting cases, the two definitions coincide.
\end{remark}

The main result of this section is the following
 \begin{theorem}\label{Bernard}~\\
Let $\mathcal D=(\widetilde D_1,...,\widetilde D_k)\in {\mathcal B_k}$
be any minimal partition associated with $\mathfrak L_k$ and let
$(\widetilde \phi_i)_i$ be any set of  positive eigenfunctions
normalized in $L^2$ corresponding to $(\lambda_1(\widetilde
D_i))_i$. Then, there exist $a_i\geq 0$, not all vanishing, such that the
functions $\widetilde u_i=a_i \widetilde \phi_i$ verify in $\Omega$
 the differential inequalities in the distributional sense
\begin{itemize}
\item[(I1)] $-\Delta \widetilde u_i\leq {\mathfrak L}_k \widetilde u_i$\,,
  $\forall i=1,\dots,k\,$,
\item[(I2)] $-\Delta\left( \widetilde u_i-\sum_{j\neq i}\widetilde
    u_j\right)\geq {\mathfrak L}_k\left(\widetilde u_i-\sum_{j\neq i}
    \widetilde u_j\right)$,
  $\forall i=1,\dots,k\,$.
\end{itemize}
\end{theorem}
\begin{remark}~\\
Note that at this stage we do not know if the  $\widetilde D_i$'s
 are connected and consequently if the $\widetilde \phi_i$'s
 are unique.  It will be shown in the next section that these
 properties are true in
two dimensions.  
\end{remark}
 
The following results were proved in \cite{CTV:2005}:
\begin{theorem}\label{extremality_condition_p}~\\
Let $p\in[1,+\infty)$ and let $\mathcal D=(D_1,...,D_k)\in {\mathcal B_k}$ be a minimal partition associated with $\mathfrak L_{k,p}$ and let $(\phi_i)_i$ be any set of  positive eigenfunctions
normalized in $L^2$ corresponding to $(\lambda_1(D_i))_i$. Then, there exist $a_i>0$, such that the
functions $u_i=a_i\phi_i$ satisfy in $\Omega$ the differential
inequalities
 in the distribution sense
\begin{itemize}
\item[(I1)] $-\Delta u_i\leq \lambda_1(D_i)u_i$,
\item[(I2)] $-\Delta\left( u_i-\sum_{j\neq i}u_j\right)\geq \lambda_1(D_i)u_i-\sum_{j\neq i} \lambda_1(D_j)u_j$.
\end{itemize}
\end{theorem}
\begin{remark}\label{osservazioni}~\\
In particular, this implies that $U=(u_1,...,u_k)$ is in the class
$\mathcal S^*$ as defined in \cite{CTV2}. Hence  Theorem~8.3 in \cite{CTV2}
 ensures the Lipschitz continuity of the $u_i$'s in the interior of
 $\Omega$. 
Therefore we can choose a partition made of open
 representatives $D_i=\{u_i>0\}$.  
\end{remark}
Moreover, taking the limit as $p\to +\infty$, the following result
 was shown in \cite{CTV:2005}:
\begin{theorem}\label{limit}~\\
There holds
$$
\lim_{p\to +\infty}\mathfrak L_{k,p}= \mathfrak L_k\,.
$$
 Moreover, there exists a minimizer of $\mathfrak L_k$ such that (I1)-
 (I2) hold for suitable non negative multiples $u_i=a_i\phi_i$ of an appropriate set of
 associated eigenfunctions.
\end{theorem}

{\bf Let us start the proof of Theorem \ref{Bernard}.}\\

Let $(\widetilde D_1,...,\widetilde D_k)\in {\mathcal B_k}$ be a
 particular minimal partition 
 associated with $\mathfrak L_k$ and let
 $(\widetilde \phi_1,...,\widetilde\phi_k)$ be any choice of associated
 eigenfunctions. We wish to prove
 that (I1)-(I2)
 hold for a suitable set of multiples of the $\widetilde \phi_j$'s.
 We consider, for a given 
\begin{equation}\label{condq}
q\in (1,N/(N-2))\;,
\end{equation}
 (or $q\in (1,+\infty)$ when $N=2$),
 the penalized Rayleigh
 quotient:
$$
\mathcal F_{k,p}(u_1,...,u_k)=
\left(\frac{1}{k}\sum_{i=1}^k \left(\frac{\int_\Omega |\nabla
u_i(x)|^2dx}{\int_\O|u_i(x)|^2dx}\right)^p\right)^{1/p}+\sum_{i=1}^k\left(1-\frac{\int_\O
  u_i(x)^q \widetilde\phi_i(x)^q\,dx}{\left(\int_\O u_i(x)^{2q}\,
  dx\int_\O \widetilde \phi_i(x)^{2q}\, dx\right)^{1/2}}\right)\;.
$$
We consider the minimization problem 
\begin{equation}\label{minM}
\mathcal M_{k,p}=\inf\left\{\mathcal F_{k,p}(u_1,...,u_k)\;:\; (u_1,...,u_k)\in\mathcal U\right\}\;,
\end{equation}
where
\begin{equation}
\mathcal U=\left\{(u_1,...,u_k)\in (W^1_0(\O))^k\;:\;u_i\cdot
u_j=0\;,\mbox{\,for\,} i\neq j\;,u_i\geq 0\;,u_i\not\equiv 0\;,\;\forall i=1,...,k\right\}\;.
\end{equation}

We note that the condition on $q$ permits  to have (weak and strong)
 continuity and differentiability 
in $W^1_0(\Omega)$ of the penalization term, 
which involves integrals of
powers of $u_i$. This will be  used later 
 to apply the direct method 
of the Calculus of Variations 
and to differentiate 
$\mathcal F_{k,p}$ at
 the minimum.\\
 It is also worthwhile noticing that $\mathcal F_{k,p}$ is invariant
 by multiplication:
\begin{equation}\label{multinv}
\mathcal F_{k,p}(a_1u_1,...,a_ku_k)=\mathcal F_{k,p}(u_1,...,u_k)\;, \forall a_i\neq 0\;.
\end{equation} 

Recalling Definition \ref{eigenvalue} we have:
\begin{proposition}\label{levels}~\\
 There holds, for every $p\in[1,+\infty)$,
$$\frac{1}{k^{1/p}}\mathfrak L_k\leq \mathfrak L_{k,p}
\leq\mathcal M_{k,p}\leq \mathfrak L_k\;.
$$
\end{proposition}
{\bf Proof}~\\
It is an immediate consequence of Jensen and H\"older inequalities. 

\begin{lemma}~\\ For every $p\in[1,+\infty)$,  the value $\mathcal M_{k,p}$ is achieved.
\end{lemma}
{\bf Proof} ~\\
 Using the invariance by multiplication \eqref{multinv},
 we can choose a bounded minimizing sequence, having as weak limit the configuration $(u_1,...,u_k)\in\mathcal U$.
Now the assertion simply follows from the weak lower semi-continuity of the norm and the compact embeddings of $W^1_0(\O)$ into $L^s(\O)$ 
for any $s\in[1,+\infty)$, whenever $N=2$, and for any $s \in[1,2N/(N-2))$ when $N\geq 3$.

\begin{lemma}\label{keylemma}~\\
Let $\Lambda >0$ and 
let $U=(u_1,...,u_k)$ be any minimizer of $\mathcal M_{k,p}$ normalized in such a way that
\neweq\label{normalization}
\left(\int_\O|\nabla u_i|^2\,dx\right)^{p-1}=
\left(\Lambda\int_\O|u_i|^2\, dx\right)^{p}\;,\qquad \forall i=1,...,k\;.
\endeq
Define
\neweq\label{fi}
f_i(u)(x)=\frac{-\gamma q}{2\left(\int_\O u(x)^{2q}dx\int_\O \widetilde \phi_i(x)^{2q}dx\right)^{1/2}}
\left[u(x)^{q-1}\widetilde\phi_i(x)^q-\frac{\int_\O u(x)^{q}\widetilde\phi_i(x)^qdx}{\int_\O u(x)^{2q}dx}
u(x)^{2q-1}\right]\;,
\endeq
where
\neweq\label{gamma}
\gamma=\Lambda^{-p}\left(\frac{1}{k}\sum_{i=1}^k \left(\frac{\int_\Omega |\nabla
u_i(x)|^2dx}{\int_\O|u_i(x)|^2dx}\right)^p\right)^{1-1/p}\;.
\endeq

Then $U$ satisfies the differential inequalities in the distribution sense
\begin{enumerate}
\item[(I1)] $-\Delta u_i\leq \lambda_1(D_i)u_i+ f_i(u_i)$,
\item[(I2)] $-\Delta\left( u_i-\sum_{j\neq i}u_j\right)\geq \lambda_1(D_i)u_i+ f_i(u_i)-\sum_{j\neq i} \left(\lambda_1(D_j)u_j+ f_j(u_j)\right)$.
\end{enumerate}
\end{lemma}

{\bf Proof}~\\
 For a fixed index $i$, 
let us introduce
$$
\widehat u_i=u_i-\sum_{j\neq i}u_j\;.
$$
Let  $\varphi\geq 0$, $\varphi\in W^1_0(\Omega)$, 
and, for $t>0$ very small, let us 
  define a new test function
$V=(v_1,\dots,v_k)$, belonging to $(W_0^1(\Omega))^k$, as follows: 
$$ 
v_j=\left\{
\begin{array}{ll}
\left(\widehat u_i+t\varphi\right)^+, &\mbox{if $j=i$}\,,\\
\left(-u_j+t\varphi\right)^-=
\left(\widehat u_i+t\varphi\right)^-\chi_{u_j>0}, &\mbox{if $j\neq i$}\,.
\end{array}\right.
$$
We first remark that there is differentiability (with respect to $t$) of  all the terms
which do not involve derivatives. Indeed, since the map $u\to (u^+)^r$
is differentiable,
 we have, for any set of functions $\eta_j\in L^s(\Omega)$ and $r>1$:
$$\int_\O \eta_j v_j^r\,dx =\left\{
\begin{array}{ll}
\int_\O \eta_j u_j^r\,dx\;+rt\int_\O \eta_j u_j^{r-1}\varphi\, dx\;+o(t),
&\mbox{if $j=i$}\,,\\
\int_\O \eta_j u_j^r\,dx\;-rt\int_\O \eta_j u_j^{r-1}\varphi\,dx \; +o(t),
&\mbox{if $j\neq i$}\,.\\
\end{array}\right.
$$
By the Sobolev Embedding Theorem, this expansion holds with respect to the $W^1_0(\Omega)$--norm provided $s\in(1,+\infty]$
and $r\leq (1-1/s)(2N/(N-2))$.  As a first application, letting
$$\alpha_j=\frac{1}{t}\left\{\int_\O|v_j|^2\,dx\;-\int_\O| u_j|^2\,dx\right\}\;,$$
 and $r=2$, 
we have
$$\alpha_j=\left\{
\begin{array}{ll}
2\int_\O u_j\varphi\,dx\;+o(1),&\mbox{if $j=i$}\,,\\
-2\int_\O u_j\varphi\,dx\;+o(1),&\mbox{if $j\neq i$}\,.\\
\end{array}\right.
$$
Moreover, letting
$$\beta_j=\frac{1}{t}\left\{\left(1-\frac{\int_\O v_j(x)^q\widetilde\phi_j(x)^qdx}{\left(\int_\O v_j(x)^{2q}dx\int_\O \widetilde \phi_j(x)^{2q}dx\right)^{1/2}}\right)-\left(1-\frac{\int_\O u_j(x)^q\widetilde\phi_j(x)^qdx}{\left(\int_\O u_j(x)^{2q}dx\int_\O \widetilde \phi_j(x)^{2q}dx\right)^{1/2}}\right)\right\}\;,$$
we find, recalling that $q\in(1,N/(N-2))$,
 by the usual differentiation rules
$$\beta_j=\left\{
\begin{array}{l}
\frac{-q}{\left(\int_\O u_j(x)^{2q}dx\int_\O \widetilde \phi_j(x)^{2q}dx\right)^{1/2}}\left[\int_\O u_j(x)^{q-1}\widetilde\phi_j(x)^q\varphi(x)dx-\frac{\int_\O u_j(x)^{2q-1}\varphi(x)dx\int_\O u_j(x)^{q}\widetilde\phi_j(x)^qdx}{\left(\int_\O u_j(x)^{2q}dx\right)}\right]\\+o(1),\quad\mbox{if $j=i$},\\
\\
\frac{+q}{\left(\int_\O u_j(x)^{2q}dx\int_\O \widetilde \phi_j(x)^{2q}dx\right)^{1/2}}\left[\int_\O u_j(x)^{q-1}\widetilde\phi_j(x)^q\varphi(x)dx-\frac{\int_\O u_j(x)^{2q-1}\varphi(x)dx\int_\O u_j(x)^{q}\widetilde\phi_j(x)^qdx}{\left(\int_\O u_j(x)^{2q}dx\right)}\right]\\+o(1),\quad\mbox{if $j\neq i$}\,.\\
\end{array}\right.
$$
On the other hand, differentiation with respect to $t$ may fail when we consider the gradient integrals. Let us denote

$$\delta_j=\frac{1}{t}\left\{\int_\O|\nabla v_j|^2\,dx\;
-\int_\O|\nabla u_j|^2\,dx\right\}\;.$$

Although $t\delta_j\to 0$ as $t\to 0$, the $\delta_j$'s themselves
 can be unbounded, in general, for they involve boundary integrals 
which are not necessarily finite for functions in  $W^1_0(\Omega)$. On the other hand, from the definition

$$
\int_\Omega|\nabla v_j|^2\,dx\;
\leq \int_\Omega|\nabla (u_j-t\varphi)|^2\,dx  \;,\;\quad\mbox{if $j\neq i$,}$$
 
we can easily deduce  that

\neweq\label{deltasneg}
\delta_j\leq -2\int_\O\nabla u_j\cdot \nabla\varphi\,dx\; +o(1)\;,
\;\quad\mbox{if $j\neq i$, and $\varphi\geq 0$}\;,
\endeq 
while, from
 $$
\begin{array}{ll}
t\sum_j \delta_j &=\sum_j\left\{\int_\O|\nabla v_j|^2\,dx\;
-\int_\O|\nabla u_j|^2\,dx \right\}\\
&
=\int_\O|\nabla (\widehat u_i+t\varphi)|^2\,dx\;
-\int_{\widehat u_i+t\varphi=0}|\nabla\widehat u_i+t\varphi|^2\,dx\;
-\int_\O|\nabla \widehat u_i|^2\,dx\;
+\int_{\widehat u_i=0}|\nabla\widehat u_i|^2\,dx\\
& 
=2t\int_\O\nabla\widehat u_i\cdot \nabla \varphi\,dx \;
+ t^2\int_\O|\nabla\varphi|^2\, dx\;,
\end{array}$$
we easily conclude that
\neweq\label{sommadelta}
\sum_j \delta_j=2\int_\O\nabla\widehat u_i\cdot \nabla \varphi\,dx\;
 +o(1)\;.
\endeq

Let us estimate, for a fixed index $j$, the difference:
 $$
\left(\frac{\int_\Omega |\nabla
v_j(x)|^2dx}{\int_\O|v_j(x)|^2\,dx}\right)^p-
\left(\frac{\int_\Omega |\nabla
u_j(x)|^2\,dx}{\int_\O|u_j(x)|^2\,dx}\right)^p=
pt\Lambda^p\left(\delta_j-\lambda_1(D_j)\alpha_j+o(\delta_j)\right)\;,$$
 
here we used the normalization condition \eq{normalization}, which implies
\neweq\label{normalization2}
\int_\O|\nabla u_j|^2\,dx =\left(\frac{\lambda_1( D_j)}{\Lambda}\right)^{p}\;.
\endeq

 On the other hand, we have:
 $$\left(1-\frac{\int_\O v_j(x)^q
\widetilde\phi_j(x)^q\,dx}{\left(\int_\O v_j(x)^{2q}\,dx\;
\int_\O \widetilde \phi_j(x)^{2q}\,dx\right)^{1/2}}\right)
-\left(1-\frac{\int_\O u_j(x)^q \widetilde\phi_j(x)^q\, dx}
{\left(\int_\O u_j(x)^{2q}\,dx\int_\O \widetilde \phi_j(x)^{2q}\,dx\right)^{1/2}}\right)=t\beta_j\;.$$


Now we prove inequality (I1). We select $j\neq i$ and we replace only
 the $j$'th component $u_j$ with $v_j$. 
 We obtain, as $t\to 0^+$,

$$
\begin{array}{ll}
        0&\leq \frac{1}{t}\left(\mathcal F_{k,p}(u_1,...,v_j,...,u_k)-\mathcal F_{k,p}(u_1,...,u_j,...,u_k)\right)\\
        &=\frac{1}{\gamma}\left(\delta_j-\lambda_1(D_j)
\alpha_j-\gamma\beta_j\right)+o(\delta_j)+o(1)\\
\end{array}
$$

This inequality, and the boundedness of the $\alpha_j$'s, the
$\beta_j$'s 
and $\gamma$ gives a lower bound  of the $\delta_j$'s. On the other hand \eq{deltasneg} gives an upper bound  of the $\delta_j$'s, which are consequently 
 bounded as $t\to 0$. 
Hence we can
deduce from \eq{deltasneg}  and the last inequality that
$$
0\leq -2\int_\O\left(\nabla u_j\cdot \nabla\varphi
  -\lambda_1( D_j)u_j\varphi-f_j(u_j)\varphi\right)\,dx\;.
$$

Since this holds for every pair of indices $i\neq j$ (though here $i$
does not appear) and every non negative test function $\varphi$, inequality (I1) is proved.

To prove inequality (I2), we argue by contradiction and we assume the existence of $\varphi\geq 0$ such that
$$\int_\O\nabla \widehat u_i\cdot \nabla\varphi\,dx<
\int_\O\left(\lambda_1(D_i)u_i(x)+ f_i(u_i)(x)-\sum_{j\neq
    i}\lambda_1( D_j)u_j(x)+f_j(u_j)(x)\right)\varphi(x)\,dx\;,$$
or, in other words, 
\neweq\label{testineq}
 \int_\O\nabla \widehat u_i\cdot \nabla\varphi\,dx 
<\sum_{i}\left(\lambda_1( D_j)\frac{\alpha_j}{2}+\gamma\frac{\beta_j}{2}\right)+o(1)\;.
 \endeq
  

Now, by the minimization property of $U$, we have,
$$
0\leq \mathcal F_{k,p}(v_1,...,v_k)-\mathcal
F_{k,p}(u_1,...,u_k)=\frac{t}{\gamma}\sum_j\left(\delta_j-\lambda_1( D_j)\alpha_j-\gamma\beta_j\right)+o(t)\;,
$$
in contradiction with \eq{testineq} and \eq{sommadelta}.

Theorem \ref{Bernard} will easily follow from the next two results:
\begin{lemma}\label{convergence}~\\
As $p\to +\infty$ any family of minimizers 
of \eqref{minM} satisfying the normalization condition
 \eq{normalization} with $\Lambda=\mathfrak L_{k}$
 converges, up to a subsequence, to a multiple 
$(a_1\widetilde\phi_1,...,a_k\widetilde\phi_k)$
 ($a_i\geq0$, not all vanishing) strongly in
$(W^1_0(\Omega))^k$. 
\end{lemma}

{\bf Proof}~\\
>From Proposition \ref{levels} we have, for the minimizers $u_{i,p}$
and the corresponding $D_{i,p}$ (we now mention the reference to $p$
 which will then tend to $+\infty$),
$$
\frac{1}{k}\sum_i\int_\O|\nabla
 u_{i,p}|^2\,dx=\frac{1}{k}
\sum_i\left(\frac{\lambda_1(D_{i,p})}{\Lambda}\right)^p\leq 
 \left(\frac{\mathcal M_{k,p}}{\mathfrak  L_k}\right)^p\leq 1\,,
$$
while
$$
\frac{1}{k}\sum_i\int_\O|\nabla u_{i,p}|^2\,dx\geq
\left(\frac{\mathfrak L_{k,p}}{\mathfrak  L_k}\right)^p\geq
\frac{1}{k}\;.
$$

Hence the family is bounded in $(W^1_0(\Omega))^k$ and does not vanish.  We
extract a sequence $(u_{i,p_n})_{n\in \mathbb N}$ possessing a limit,
in the weak $(W^1_0(\Omega))^k$--topology and in any $(L^r(\Omega))^k$, for
subcritical $r$'s. We denote
 by  $(\widetilde u_i)_{i=1,\dots,k}$ this limit.
We infer that the weak limit cannot be identically zero. We have indeed:
$$
\lambda_1(D_{i,p_n})=\frac{\int_\O|\nabla u_{i,p_n}|^2 \,
  dx}{\int_\O|u_{i,p_n}|^2\, dx }\leq
k^{1/p_n}\mathfrak L_k\;,\;\forall i=1,...,k\;,
$$
so that
$$\frac{1}{k}\sum_{i}\int_\O|u_{i,p_n}|^2\,dx\geq
\frac{1}{k^{1/p_n}\mathfrak  L_k}
\sum_{i}\int_\O|\nabla u_{i,p_n}|^2\, dx\;.$$
We further remark that, if for some $i$ the weak limit happens to be zero, then the strong limit vanishes too. \\
We claim that, for suitable non negative $a_i$'s we have $\widetilde
u_i=a_i\widetilde \phi_i$. This is obvious if $\widetilde u_i\equiv 0$. If not, since 
by Proposition~\ref{levels} 
 $$
\lim_{p\to +\infty}\left( \mathcal M_{k,p}-\mathfrak  L_{k,p}\right)=0\;,
$$
 we deduce that, whenever $\widetilde u_i\not\equiv 0$, 
$$
1-\frac{\int_\O \widetilde u_i(x)^q
  \widetilde\phi_i(x)^q\,dx}{\left(\int_\O \widetilde u_i(x)^{2q}\,dx\int_\O
  \widetilde \phi_i(x)^{2q}\,dx\right)^{1/2}}=0\;,
$$
and therefore that $\widetilde u_i$ is a multiple of $\widetilde
\phi_i$. To pass from weak to strong convergence we first notice that
each $u_{i,p_n}-\widetilde u_i$ converges weakly and strongly in $L^2$
to zero. 
Now we recall that $u_{i,p_n}$ satisfy inequalities (I1)--(I2) of Lemma~\ref{keylemma}.
Let us multiply (I1) 
  by $(u_{i,p_n}-\widetilde u_i)^+$,
 (I2) by $(u_{i,p_n}-\widetilde u_i)^-$ and take the
 difference. We obtain
$$
\int_\Omega\nabla u_{i,p_n}\cdot \nabla (u_{i,p_n}-\widetilde
u_{i})\,dx 
\leq \int_{\Omega}\sum_{j\neq i}\nabla u_{j,p_n}\cdot \nabla
(u_{i,p_n}-\widetilde u_i)^-\,dx\; + o(1)\;.
$$
Since $u_{i,p_n}(x)u_{j,p_n}(x)$ vanishes almost everywhere
 and $\widetilde u_i\geq 0$, we infer
 $$
\begin{array}{ll}
\int_\Omega\nabla u_{i,p_n}\cdot \nabla (u_{i,p_n}-\widetilde u_i)\, dx
&\leq
\sum_{j\neq i}\int_{\Omega}\nabla u_{j,p_n}\cdot \nabla
\widetilde u_i\, dx\;  +o(1)\\& =\sum_{j\neq i}\int_{\Omega}\nabla \widetilde
u_j\cdot \nabla
\widetilde u_i\, dx \;+o(1)=o(1)\;.
\end{array}
$$

Thus we can deduce strong convergence from the weak. We have  indeed  
$$
|| \nabla ( u_{i,p_n}-\widetilde u_i)||^2
=\int_\Omega\nabla u_{i,p_n}\cdot \nabla (u_{i,p_n}-\widetilde u_i)\,dx
 - \int_\Omega\nabla \widetilde u_{i}\cdot \nabla (u_{i,p_n}-\widetilde u_i)\,dx
 =  o(1)\;.
$$
 \begin{lemma}~\\
Let $U_n=(u_{1,p_n},...,u_{k,p_n})$ ($n\in \mathbb N$) as in the proof of Lemma~\ref{convergence}. Then its limit, as $n\rightarrow +\infty$,  
 $\widetilde{U} :=
(\widetilde u_1,\dots,\widetilde u_k)$ verifies the inequalities in the statement of Theorem~\ref{Bernard}. 
\end{lemma}
{\bf Proof}~\\
First of all, we wish to pass to the limit in formulas (I1)-(I2) of
Lemma~\ref{keylemma}, in the sense of distributions. From the previous
lemma we deduce that $-\Delta u_{i,p_n}\to -\Delta \widetilde u_i$  in
$H^{-1}(\Omega)$\,. 
Hence one can pass to the limit in inequality (I1). Let us turn to (I2). We remark that 
also $f_i(u_{i,p_n})$
converge to $f_i(\widetilde u_i)$, provided $\widetilde u_i\not\equiv
0$.
 For such $i$'s inequality (I2) passes to the limit, 
because so does its right hand. On the other hand, 
when the limit $\widetilde u_i$ does vanish then (I2) holds because 
of (I1) and the fact that $-\Delta\widetilde u_i=0$. 
In order to end the proof, we have to prove convergence
 of the eigenvalues $\lambda_1(D_{i,p_n})$ to $\mathfrak  L_k$ 
whenever $\widetilde u_i$ does not vanish identically. At first
 we notice that the $\lambda_1(D_{i,p_n})$'s do converge,
 thanks to the strong convergence of the $u_{i,p_n}$'s
 to limits $\lambda_1(\widetilde D_i)\leq {\mathfrak  L}_k$. 
Assuming $\lambda_1(\widetilde D_i)< {\mathfrak  L}_k$ 
we deduce from \eq{normalization2}, 
using again the strong convergence, that $\widetilde u_i \equiv 0$.

\begin{remark}\label{osservazionia}~
\begin{enumerate}
\item[(a)] Thanks to \cite{CTV2}, Theorem 8.3
 all the $u_{i,p}$'s and their limits  $\widetilde u_i$ are locally Lipschitz continuous in the interior of $\O$  
 and continuous up to the boundary, for (I1), if the boundary
 $\partial\Omega$ is Lipschitz, or has the interior 
 cone property; moreover, they are globally Lipschitz up to the
 boundary, if $\partial\Omega$ is of class $\mathcal C^1$.
\item[(b)] Of course, since the eigenfunctions are normalized in $L^2$
$$a_i\mathop{=}^{def}\Vert \widetilde u_i\Vert_2\,.$$ 
In general, it may happen that some of the $a_i$'s vanish; we denote
\begin{equation}\label{knot}
{\bf k_0}=\{i\in\{1,...,k\}: a_i= 0\}\,.
\end{equation}
\item[(c)] Going back to the proof of Lemma~\ref{convergence}, we can
  extract a subsequence with the further property that  the $u_{i,p_n}/\Vert
  u_{i,p_n}\Vert_2$ converge strongly in $L^2(\Omega)$ and weakly in 
$W^1_0(\Omega)$ to
  $\widetilde\phi_i$, also for those indices $i$'s for which the
  component $u_{i,p_n}$ normalized as in \eq{normalization2} strongly
  converges to $ 0$.
\item[(d)]  We also infer from \eq{normalization2} that
$$\lambda_1(\widetilde D_i)\leq\mathfrak L_k\;\qquad\forall i\in\{1,...,k\}\;,
$$
while
$$\lambda_1(\widetilde D_i)=\mathfrak  L_k\;,\qquad\forall i\not\in{\bf k_0}\;.
$$
\item[(e)] Obviously the system of differential inequalities (I1)-(I2)
  are fulfilled by the set $(\widetilde u_i)$ with $i\not\in \bf k_0$\;.
\item[(f)]  For $i\not\in \bf k_0$, the $\widetilde D_i$'s are open and
  possess finitely many connected components.

\end{enumerate}

\end{remark}
\section{Further results in two dimensions}\label{susanna2}
In this section, we recall from the introduction that
 we work under  Assumptions
 \ref{Asswhole1} and \ref{Asswhole2}. 
This  implies that $\partial\O$ has 
finitely many connected components. For 
 simplicity, we shall omit the potential $V$  in the
following discussion.  All the arguments can be straightforwardly
extended (sometimes at the price to replace $\mathcal C^2$ by
 $\mathcal C^{1,1-}:=\cup_{\alpha<1} \mathcal C^{1,\alpha}$) in
order to cover the case of a non vanishing  bounded potential. A special
caution is only due in the proof of Theorem \ref{grad.nullo},    where the needed
 extra argument is outlined . 

\subsection{Case when $
\mathbf k_0=\emptyset$.}
In this section we discard in a first step 
 all the identically vanishing components. Hence, from now on, 
we will assume that
$$
\mathbf k_0=\emptyset\;.
$$
Re--labeling if necessary, and taking a smaller $k$, we assume that
the components of 
 $$
U=(u_1,\dots,u_k)
$$ are non negative, non vanishing 
 $W^1_0(\Omega)$--functions, such that $u_i(x)u_j(x)\equiv 0$ almost everywhere in 
 $\O$ (for $i\neq j$),
 satisfying the two differential inequalities
 (I1)-(I2) of Theorem \ref{Bernard}. As a consequence of Remark
 \ref{osservazionia}~(a), they are continuous on the
 closure of $\O$. Hence, by expanding the set of indices 
 if necessary, we
 can always assume that the sets
$$
D_i=\{u_i>0\}$$
are  open and connected. \\
Let us define the set of
zeroes  of $U$ as
$$
{\cal Z}=\{x\in\Omega: u_i(x)=0\;,\;\forall\;i=1,...,k\}\;,
$$
and define the {\it multiplicity of $x$}$\in {\Omega}$ as the number $m(x)$:
\neweq\label{emme}
m(x)=\sharp\{i\; :\; meas(\{u_i>0\}\cap B_r(x))>0\;,\; \forall r>0\}\;.
\endeq
We shall denote by
$$
{\mathcal Z}_h=\{x\in\Omega:m(x)\geq h\}
$$
the set of points of multiplicity greater than or equal to the integer $h$ and
by
$$
{\cal Z}^h=\{x\in \Omega : m(x)=h\}\,.
$$

We remark that, by definition, $\mathcal Z^0$ is open.
Let us  now consider  $\mathcal Z^1$. This is the object of:
\begin{proposition}\label{zeta1}~\\
$$\mathcal Z^1=\bigcup_{i=1}^k D_i\;.$$
\end{proposition}
{\bf Proof}~\\
If  a ball $B_r(x_0)$ intersects only $D_i$ we deduce from (I1)--(I2) that $u_i$ is a non negative solution to the differential equation $-\Delta u =\Lambda u$  on $B_r(x_0)$. By the strong maximum principle then $u_i$ is strictly positive on $B_r(x_0)$ and therefore
$B_r(x_0)\subset D_i$ and all the points of $B_r(x_0)$ have multiplicity equal to one.
\endproof
 
To continue the analysis of  the topological properties related to
points of multiplicity two and more, we shall consider the simplicial
homology groups with coefficients
  in $\mathbb Z_2$,
$H_{n}(X)$, for $n=0,1$.  
We recall that $ rank\,(H_0(X))$ is the number of connected components of
$X$. For open subsets of an euclidean space, as the fundamental group $\pi_1$ is
 already abelian, there holds $rank(\pi_1(X))=rank (H_1(X))$. Finally for planar bounded open subsets,
 $rank (H_1(X))+1$ is the number of connected components of $\partial X$. A reference 
 book for  the  algebraic topology concepts is the Greenberg and Harper's book \cite{gh}.

\begin{proposition}\label{finiteh1i}~\\
 If $H_1(\Omega)$ is finite,  so is $H_1(D_i)$, for every $i$.
\end{proposition}
{\bf Proof}~\\
Let us consider a loop $\gamma\subset D_i$  which is homotopically
trivial  
in $\Omega$ but not in $D_i$; hence denoting by $\Sigma$ the inner
region of $\gamma$, 
we have that $\Sigma\subset\Omega$ but $\Sigma\not\subset D_i$. Let
$j\neq i$;
 then either $D_j$ is contained in $\Sigma$ or in its complement, for
 it is connected.
 To prove the proposition, we argue by contradiction and we assume
 that $H_1(D_i)$ is infinite.
 Then also $rank(\pi_1(D_i))$ is infinite; thus we infer the existence 
of at least one loop $\gamma=\partial \Sigma$ such that
 $\Sigma\cap \bigcup_{j\neq i}D_j=\emptyset$ and $\Sigma\not\subset D_i$.
Since all the $u_j$'s ($j\neq i$) vanish identically in $\Sigma$, we
deduce from (I1)-(I2) that $- \Delta u_i=\Lambda u_i$ in $\Sigma$ and,
by the strong maximum principle,  that $u_i$ is strictly positive
there;
 thus $\Sigma\subset D_i$, a contradiction.\endproof

Let us consider
$$
\Gamma_{i,j}=\partial D_i\cap\partial D_j\cap {\cal Z}^2\;,$$
$$D_{i,j}=D_i\cup D_j\cup\Gamma_{i,j}\;.$$
Our next goal consists in showing that the $\Gamma_{i,j}$ consist in a finite
number of (possibly open) arcs. This will require some topological considerations.

\begin{proposition}\label{localc1}~\\
Let $x_0\in \Omega$ such that $m(x_0)=2$. Then $u_i-u_j$ is in $\mathcal C^{1,1-}$ in some neighborhood of $x_0$
 and $\nabla (u_i-u_j)(x_0)\neq
0$. Furthermore ${\cal Z}^2$ is locally a $\mathcal C^{1,+}$--curve through $x_0$.
\end{proposition}
{\bf Proof}~\\
 ~Relabeling, we can always assume 
 $x_0\in \partial\{u_1>0\}\cap\partial\{u_2>0\}$; thus
for all $r$ small enough  $B(x_0,r)\cap D_i=\emptyset$ for
all $i>2$. Then $u=u_1-u_2$ satisfies the equation
$-\Delta u=\Lambda u$ in $B(x_0,r)$ 
and is consequently $\mathcal C^{1,1-}$ near $x_0$. Therefore  the
 zero set of $u$ near $x_0$ (by a standard result
 on the zero set) is made up by a finite number of even 
regular curves starting from $x_0$. But there are actually
 only two arcs meeting at $x_0$. Indeed, if not,
 at least one of the $D_i$'s should be 
disconnected. In this case the zero set is actually locally a regular line
 passing through $x_0$ and the Boundary Point  Lemma gives the proposition.
\endproof

\begin{proposition}\label{finiteh1ij}~\\
If $H_1(\Omega)$ is finite so is $H_1(D_{i,j})$, for every $i\neq j$.
\end{proposition}
{\bf Proof}~\\ 
This is an obvious statement if $\Gamma_{i,j}=\emptyset$. 

Let us consider a loop $\gamma\subset D_{i,j}$ which is trivial in $\Omega$
 but not in $D_{i,j}$; hence $\gamma=\partial \Sigma$  
with $\Sigma\subset\Omega$ but $\Sigma\not\subset D_{i,j}$.
 Let $\ell \not\in \{i,j\}$; then either $D_\ell $ is contained in $\Sigma$
 or in its complement. Thus, arguing by contradiction we infer
 the existence of a loop $\gamma=\partial \Sigma$ such that
 $\Sigma\cap \bigcup_{\ell \neq i,j}D_\ell =\emptyset$ and
 $\Sigma\not\subset D_{i,j}$. Since all the $u_\ell$'s ($\ell \neq i,j$) 
vanish identically in $\Sigma$, we deduce from (I1)-(I2) that 
$-\Delta (u_i-u_j)=\Lambda (u_i-u_j)$ in $\Sigma$ and, 
from Proposition \ref{localc1} and the unique continuation principle
 we infer that $\mathcal Z\cap \Sigma\subset \Gamma_{i,j}$;
 thus $\Sigma\subset D_{i,j}$, a contradiction.\endproof

\begin{lemma}~\\
 If $H_1(D_{i,j})$ is finite, then $\Gamma_{i,j}$ has finitely many connected components.
\end{lemma}
{\bf Proof}~\\
To prove the statement we will take advantage of the
Mayer--Vietoris Theorem. The Mayer--Vietoris sequence is usually  proven to be exact for a triad $X$, $X_1$ and $X_2$ where 
$X_1$ and $X_2$ are open subsets of the topological space $X$ and $X$ is the 
union of $X_1$ and $X_2$ (such a triplet is called and admissible triad). 
Here we would like to apply the Mayer--Vietoris sequence to the triad
$X=D_{i,j}$, $X_1=D_i\cup\Gamma_{i,j}$ and $X_2=D_j\cup\Gamma_{i,j}$, but these
latter two are not open in $D_{i,j}$. However, the Mayer--Vietoris sequence is still 
available because, thanks to Proposition \ref{localc1}, the $\Gamma_{i,j}$'s 
are regular embedded one--dimensional submanifolds in $D_{i,j}$. Hence
each $D_\ell \cup\Gamma_{i,j}$ 
is a euclidean neighborhood retract in $D_{i,j}$ and, as such it has the same homology as 
the corresponding $D_\ell $ ($\ell =i,j$). 
Thus, following \cite{gh}, the triplet $D_{i,j}$,  $D_i\cup\Gamma_{i,j}$, $D_j\cup\Gamma_{i,j}$ is a proper
excision triad  and thus the Mayer--Vietoris  sequence is exact:
$$
H_1(D_{i,j})\stackrel{\partial_*}{\rightarrow}H_0(\Gamma_{i,j})
\stackrel{i_*\oplus-j_*}{\rightarrow}H_0(D_i\cup\Gamma_{i,j})
\oplus H_0(D_j\cup\Gamma_{i,j})$$
The assertion
 then follows as a consequence of Propositions
 \ref{finiteh1i}
 and \ref{finiteh1ij} 
taking into account of the connectedness of $D_\ell $ ($\ell =i,j$). Indeed since both $H_1(D_{i,j})$ and $H_0(D_i\cup\Gamma_{i,j})
\oplus H_0(D_j\cup\Gamma_{i,j})=\mathbb Z_2\oplus \mathbb Z_2$ are
 finite, 
we have that both the range of ${\partial_*}$
and the coker of ${i_*\oplus-j_*}$ are finite and thus so is $H_0(\Gamma_{i,j})$. \endproof

The following results follow from  \cite{CTV2}:

\begin{theorem}\label{grad.nullo}~\\
We have, $\forall i=1,\dots,k$~:
 \begin{itemize}
\item[(a)] $u_i\in W^{1,\infty}_{loc}(\Omega)$.
\item[(b)] The function $\sum_{i=1}^k\vert \nabla u_i\vert$ admits a continuous representative in $\Omega$.
\item[(c)] If $x_0\in{\cal Z}_3$ then, 
$$
\lim_{{x\to x_0}}\sum_{i=1}^k\vert \nabla u_i(x) \vert=0\;.
$$
\item[(d)] If $x_0\in{\cal Z}_3$ then,
$$
\lim_{x\to x_0}\frac{u_i(x)}{\vert x-x_0\vert}=0\;.
$$
\end{itemize}
\end{theorem}

{\bf Proof}~\\
Part (a) is indeed stated as Theorem 8.3 of \cite{CTV2}. To prove part
(b), we first observe that locally in $\bigcup_{i=0}^2 \mathcal{Z}^i$ (which is open in $\Omega$)  the $u_i$'s satisfy the differential equation
$-\Delta(u_i-u_j)=\Lambda (u_i-u_j)$: hence they are regular (of class
$\mathcal C^1\left(\bigcup_{i=0}^2 \mathcal{Z}^i \right)$). The
continuity up to the $\mathcal{Z}_3$ (which is indeed the boundary of
$\bigcup_{i=0}^2 \mathcal{Z}^i $) is then a consequence of
the vanishing of the limit stated in part (c), which is indeed Theorem
9.3
 in \cite{CTV2}. The proof of Theorem 9.3 was originally performed
 in the absence of the $L^\infty$ potential, but 
 all the arguments can be promptly adapted to cover also this case.
 A special care is needed when, at the beginning, the function
 $\sum_{i=1}^k\vert \nabla u_i\vert^2$ is shown to be a subsolution
 to a linear differential equation. Subsequently the mean
 value property is applied. This is not exactly true in the 
present situation, for the linear problem is now perturbed 
by a term of the form $\nabla w\cdot \nabla((\Lambda+V(x))w)$,
 where $w$ is an auxiliary function (which can be taken, 
by the way, the same as defined in \eq{eq:auxiliary}). 
After integrating by parts one easily sees that the contribution 
of this term is negligible. This observation permits to follow,
 from then on, exactly the same arguments of the quoted paper.\\ 

\begin{lemma}\label{aderenza}~\\ Let $x_0\in{\cal Z}_3$.
Then there exists a sequence  $\{x_n\}\subset\Omega$ such that $m(x_n)=2$ and
$x_n\to x_0$ as $n\rightarrow +\infty$. \end{lemma} 
{\bf Proof}~\\
 Assume not, then there would be an
element $y_0$ of ${\cal Z}_3$ having a positive distance $d$ from
${\cal Z}^2$. Let $r<d/2$; then the ball $B(y_0,r)$ intersects at
least three of the $D_i$'s; therefore there exist $i\in \{1,\dots,k\}$,  $x\in D_i$
and $z_0\in{\cal Z}_3$ such that $\rho=d(x,z_0)=d(x,{\cal
Z}_3)<d(x,{\cal Z}^2)$. Then the ball $B(x,\rho)$ is tangent from
the interior of $ D_i$ to ${\cal Z}_3$ in $z_0$. Furthermore
$u_i$ solves an elliptic PDE (in the sense of \cite{GT:1983})
 and it is positive on $D_i$. Then we
thus infer from the Boundary Point
 Lemma (see Lemma 3.4 and formula (3.11) in \cite{GT:1983} ) that 
 \begin{equation}\label{gradgen}
\liminf_{h\to 0}\frac{u_i(z_0+h\nu)}{h}>0\;,
\end{equation}
where $\nu$ denotes the inner normal to $\partial B(x,\rho)$ at $z_0$,  in contrast with Theorem~\ref{grad.nullo}(c).\endproof

A simple but important consequence of this  discussion
  is the following result (the closure is
taken in $\overline{\Omega}$):
\begin{proposition}\label{finite}~\\
If $H_1(\Omega)$ is finite then $\overline{{\cal Z}_3}$ has a finite number of
connected components.
\end{proposition}
{\bf Proof}~\\
Indeed, we have $\overline{{\cal Z}_3}\subset \cup_{i,j}\left(\overline{\Gamma_{i,j}}\setminus\Gamma_{i,j}\right)$, and each
$\Gamma_{i,j}$ is the union of finitely many arcs, each of them being
homeomorphic to the real line (non compact case) or to  a circle. In
the first case they have  
 connected and closed $\alpha$ and $\omega$--limits.\\
We recall that the $\alpha$ and $\omega$--limits of a parametrized arc
 $\Gamma(t)$ are the
sets of the limit points as the parameter tends to $-\infty$ 
and $+\infty$ respectively: one
easily sees that bounded arcs have compact and connected $\alpha$ and
 $\omega$--limits. 
  Therefore each $\overline{\Gamma_{i,j}}\setminus\Gamma_{i,j}$ contributes with finitely many connected components by the previous proposition.
\endproof
In addition we have:
\begin{proposition}\label{finite0}~\\If $H_1(\Omega)$ is finite then
$\overline{{\cal Z}_3\cup \mathcal Z^0}$ has a finite number of
connected components.
\end{proposition}
{\bf Proof}~\\
 Indeed, by Theorem \ref{grad.nullo} and Lemma~\ref{aderenza},  the boundary
 $\partial\left(\overline{{\cal Z}_3\cup \mathcal Z^0}\right)$ is the
 same as $\partial\overline{{\cal Z}_3}$ and this last one has finitely many connected components.

Our next goal is the following
\begin{theorem} \label{isolz3.2}~\\
An isolated connected component of  $\overline{{\cal Z}_3\cup \mathcal Z^0}$ consists of a single point.
\end{theorem}
This theorem  implies straightforwardly,  
having in mind  that $\mathcal Z^0$ is open, that
\begin{corollary} \label{no_vanishing}~\\
 If $H_1(\Omega)$ is finite then
$
\mathcal Z^0=\emptyset\;.
$ 
\end{corollary}
{\bf Proof of Theorem \ref{isolz3.2}}~\\
 To prove the theorem we focus on a connected component
  $Y_0$ of $\overline{{\cal Z}_3\cup \mathcal Z^0}$ 
and we show that it is reduced to a single point. 
First we consider an open and connected neighbourhood  $\mathcal N$ 
 of $Y_0$ in $\mathbb R^2$ having a regular boundary and such that
 $\overline{\mathcal N}\cap \overline{{\cal Z}_3\cup \mathcal
 Z^0}=Y_0$. Since $Y_0$ is connected, we can choose $\mathcal N$ in such a way that 
 its boundary has exactly one or two connected components (each
 diffeomorphic to $S^1$), depending on 
 whether $Y_0$ disconnects $\mathbb R^2$ or not. This can be achieved
 by taking one connected component of a regular sublevel of a non
 negative $\mathcal C^\infty$ function having $Y_0$ as null set.
 Furthermore we can take the measure
 of $\mathcal N\setminus Y_0$ small enough that none
of the $D_i$'s is completely enclosed in $\mathcal N$.
 We may assume that $\partial \mathcal N$ intersects 
transversally the $\Gamma_{i,j}$'s and $\partial \Omega$.
 Thus, after possibly cutting--off some portions  of $\mathcal N$,
 we can assume 
that each oriented arc in  $\Gamma_{i,j}$ intersects 
$\partial \mathcal N$ exactly once. Then each $D_i\cap \mathcal N$
 has a finite number of connected components. 
 Moreover, we can manage to have  these intersections simply
 connected. 
 Indeed, assuming not, and arguing as in the proof of Proposition
 \ref{finiteh1i}, 
at least one of the other $D_j$'s should be entirely contained in
 $D_i\cap \mathcal N\setminus{Y_0}$, 
what we have excluded by taking the measure
 of $\mathcal N\setminus{Y_0}$ small enough.
 We can label the
 connected components of $D_i\cap \mathcal N$ and
 $(\mathbb R^2\setminus \O)\cap \mathcal N$ clockwise
 from $1$ to $h$, according with their intersection 
with $\partial \mathcal N$.  By taking a double covering of $ \mathcal
 N$, 
still denoted by $\mathcal N$, if necessary, we may assume that $h$ is an even integer.
Indeed, assuming $h$ to be odd,  we choose some  point $y_0\in Y_0$
 and   perform a double covering branched at $y_0$ (in complex notation
 $f(z)=(z-y_0)^2$). With some abuse of notation, we denote with the same symbols
 $D_i$, $\mathcal N$ and $Y_0$ and their pre--images,  and we denote
 $\tilde u_i(z)=u_i(f(z))$\,. Note that $$ \Delta \tilde
 u_i(z)=8|z-y_0|^2 \Delta u_i(f(z))\;.
$$
If the original $h$ was even we set $\tilde u_i=u_i$ and we define, 
in both cases, the auxiliary function 
\begin{equation}\label{eq:auxiliary}
w(x)=
\begin{cases}\sigma(x)\tilde u_i(x)\,,& \text{if $x$}\in D_i \\
0& \text{otherwise}\\
\end{cases}
\end{equation}
where  $\sigma$ is a sign assignment compatible
 with the partition of $\mathcal N$. The function 
$w$ satisfies the linear equation:

\begin{eqnarray}
-\Delta w&=\Lambda a(x)w\qquad &\text{in }\Omega\cap(\mathcal N\setminus Y_0)\label{equat}\\
w&=0 \quad\qquad\qquad&\text{on }(\partial\Omega\cap\mathcal N)\cup Y_0\;,\label{BC}
\end{eqnarray}
where $a\equiv 1$ if no double covering 
 has been performed,
 otherwise $a(x)=8|x-y_0|^2$ is the conformal factor. 
Of course $w$ is regular outside the singular component $Y_0$. Also, $w$ vanishes identically on the open set $N\cap \mathcal{Z}^0$; hence, thanks to 
Theorem \ref{grad.nullo} (c), $w$ is in $\mathcal C^1(\Omega\cap \mathcal
N)$

We claim that $w$ actually solves \eq{equat} in the whole 
of $\Omega\cap\mathcal N$. This is not a trivial fact,
for there is no a priori bound on the Hausdorff measure
 of the boundary $\partial Y_0$ of the singular set. 
In order to overcome this problem  we make the following construction:

\begin{proposition}\label{neighbs}~\\
There exists a family of  neighborhoods $\mathcal N_\delta \subset
\mathcal N$, 
decreasing to $Y_0$ and bounded by a finite number of regular arcs,
 with the property that 
\begin{equation}\label{boundary}
\lim_{\delta\to 0}\int_{\partial(\mathcal N_\delta\cap\O)}|\nabla w|\,ds\;=0\;,
\end{equation}
where $ds$ denotes the measure on the union of these regular arcs.
\end{proposition} 
 
 Postponing the proof of the proposition, we end the proof
 of Theorem~\ref{isolz3.2}. Testing  \eq{equat} with any function  $\varphi\in
 C_0^\infty(\mathcal N\cap\O)$, it follows by decomposing  the integration 
  on $\mathcal N \setminus \mathcal N_\delta$
 and on $\mathcal N_\delta$
 and then making an integration by parts
 for the second integral:
$$ 
 \begin{array}{ll}
 \vert \int_{\mathcal N\cap\Omega} &\left(\nabla w\cdot 
\nabla \varphi-\Lambda a(x)u\varphi\right) \,dx \;\vert\\
 &=\vert\int_{\mathcal N_\delta\cap\Omega} \left(\nabla w\nabla
 \varphi-\Lambda a(x)u\varphi\right)\,dx\; -
\int_{\partial(\mathcal N_\delta\cap\O)} 
 \varphi\nabla w\cdot \nu\; ds \vert \\
 &\leq  C (\sup_{\mathcal N_\delta\cap\Omega\cap\text{supp}(\varphi)}
(|w|+|\nabla w|)+\int_{\partial(\mathcal N_\delta\cap\O)}|\nabla w|\,ds)\;.
\end{array}
 $$
We have now to show that the right hand side is $o(1)$. Proposition
\ref{neighbs}   ensures  that, as $\delta\to
0$,
$$
\int_{\partial(\mathcal N_\delta\cap\O)}|\nabla w|\,ds 
= o(1).
$$
So it remains to show that, as $\delta\to 0$, there holds
$$
 \sup_{\mathcal N_\delta\cap\Omega\cap\text{supp}(\varphi)}
(|w|+|\nabla w|)
= o(1).
$$
This can be done by Theorem \ref{grad.nullo}. Indeed both $w$ and $\nabla w$ are uniformly continuous on
$\text{supp}(\varphi)$ which is compactly contained in $\Omega$, and the distance of $\mathcal N_\delta\cap \text{supp}(\varphi)$ to $Y_0$  tends to zero. Finally we recall that $w$ and its gradient vanish identically on $Y_0$.
 Hence $w$ solves \eq{equat} on the whole of $\mathcal
 N\cap\Omega$. By a classical local regularity result
 by Hartman and Wintner (\cite{hw1}, Corollary 1),  we know that  
interior critical points of solutions to such class of equations
 are isolated and have finite (local) multiplicity $m$, and
 satisfy \eqref{asfo}.

Hence we are left with the case when  $Y_0$ is contained in a connected component of the boundary.
We notice that, in this case, we do not need the double covering, for a sign assignment compatible
 with the partition of $\mathcal N$ always exists, since $\partial\O$ disconnects $\partial \mathcal N$. 
Now we need an extension of Theorem \ref{reg_nodal_int}, suitable
 to cover the case of domains  possessing a finite number of $\mathcal C^{1,+}$--corners:
namely Theorems \ref{reg_nodal} and \ref{reg_nodal_corner} in Section \ref{spns}.

In particular, these results  guarantee finiteness of critical points also at the boundary. 

{\bf Proof of Proposition \ref{neighbs}}~\\
 First of all, let us consider a triplet of adjacent domains separated
 by the two arcs $\Gamma_{i,j}$ and $\Gamma_{j,k}$. We claim 
that the distance of the arcs, relative to $D_j$, must vanish. In other words, we claim  that
 $$
\inf\{\int_0^1|\dot \gamma(s)|\,ds\;:\;\gamma(0)\in
\Gamma_{i,j},\,\gamma(s)\in D_j\,,\,\forall
s\in(0,1),\,\gamma(1)\in\Gamma_{j,k}\}=0\;.
$$

 Indeed, if not,  one could find in $D_j$ a ball, tangent to the
 boundary 
$\partial D_j\cap Y_0$ at, say, $x_0$  and having positive distance 
 from both $\Gamma_{i,j}$ and $\Gamma_{j,k}$. By the Boundary Point
 Lemma, at $x_0$ the gradient 
 of $w$ cannot vanish (in the weak sense of \eqref{gradgen}), in contradiction with the fact that $x_0\in Y_0$.

As a second remark, by integrating the equation $-\Delta u_i=\Lambda
u_i$ 
over the set  $\{u_i>\varepsilon\}$ and using  the Divergence Theorem,
we obtain,
$$ 
\lim_{\varepsilon\to 0^+}\int_{\partial \{u_i>\varepsilon\}}|\nabla
u_i|\,ds\leq C \int_\O u_i\,dx <+\infty\;.
$$

The Divergence Theorem is applicable because the level sets (according
to Hartman and Wintner's regularity result)  of the eigenfunctions of
$-\Delta +V(x)$ are compact and 
piecewise $\mathcal C^1$ when the potential $V$ is bounded. 
 Since the components of ${\partial \{u_i>\varepsilon\}}$ converge to $\Gamma_{i,j}$ as $\varepsilon\to 0$, in $\mathcal C^1$ as parametrized curved, we obtain that
$$ 
\int_{\Gamma_{i,j}}|\nabla w|\,ds<+\infty\;,\qquad\forall i,j\;.
$$
If the arc $\Gamma_{i,j}$ meets $Y_0$, then
we can choose an orientation for its parametrization $\gamma(t)$
   in a way that $\lim_{t\to+\infty} \mbox{dist} (\gamma(t),Y_0)=0$, and  we
have~:
$$\int_0^{+\infty}|\nabla w(\gamma(t))||\gamma^\prime(t)| dt <+\infty\;;$$
\[\begin{split}
&\gamma(t)\not\in Y_0\,,\,\qquad\forall t\geq 0\;;
\\
&\lim_{t\to+\infty} \mbox{dist}(\gamma(t),Y_0)=0\;.
  \end{split}
\]
Therefore we can conclude
$$
 \lim_{\delta\to 0}\int_{\Gamma_{i,j}\cap B_{\delta}(Y_0)}|\nabla
 w|\,ds=0\;,
\qquad\forall i,j\;,
$$
where we denoted $B_\delta (Y_0)=\{x\in \mathbb R^2\;,\;
 d(x,Y_0) < \delta\}$.

 Now we illustrate the construction of the boundary of the
 neighbourhood $\mathcal N_\delta$. 
We start by taking the component of $\Gamma_{i,j}\cap B_{\delta}(Y_0)$
 ending at $Y_0$;
 then we can jump from this arc to the next $\Gamma_{j,k}$, still
 remaining 
in $B_{\delta}(Y_0)\cap D_j$, following a segment of arbitrarily short length. We proceed in this way
passing from one arc to the next until we complete the loop.  \endproof

Next result straighforwardly follows from Theorems~\ref{reg_nodal} and \ref{reg_nodal_corner} applied to
the auxiliary function $w$ defined in 
\eq{eq:auxiliary}; it completes our analysis of the asymptotic expansion of the nodal set 
at multiple intersection points.

\begin{theorem}\label{asym}~\\
 Let $x_0\in \O$ with $m(x_0)=h\geq 3$
that is isolated in $\overline{\cal Z}_3$. Then there exist an integer  $n\geq h$ (the local multiplicity), $c\in
\bbbr\setminus 0$, 
$\theta_0\in (-\pi,\pi]$ such that
$$
\sum_{i=1}^h u_i(r,\theta)=c
r^{\frac{n}{2}}|\cos(\frac{n}{2}(\theta+\theta_0))|+o(r^{\frac{n}{2}})\;,
$$
as $r\to 0$, where $(r,\theta)$ denotes a system of polar coordinates around
$x_0$ and $n$ is the \emph{local multiplicity} of $x_0$.\\
 Moreover the nodal set in a neighbourhood of $x_0$ is the union of $n$ closed
arcs of class $\mathcal C^{1,+}$ meeting at $x_0$ and spanning angles
of
 of opening angle  $2\pi/n$. \\
If the boundary $\partial\Omega$ presents a $\mathcal C^{1,+}$--corner 
of amplitude $\alpha\pi$ at $x_0$ and the nodal set hits $x_0$ from inside
the corner, then there exist an integer $n$ and $R>0$ such that the component of
 $\overline{u^{-1}(\{0\})}\cap B(x_0,R)$ lying inside the corner
 is composed by $n$ $\mathcal C^{1,+}$-simple arcs which all end in $x_0$ 
and whose tangent lines at $x_0$ divide the sector
 into $n+1$ angles of equal opening angle  $\pi\alpha/(n+1)$. 
\end{theorem}

\subsection{General case.}
Let us come back to the general case. We no more assume a priori that 
 ${\bf k}_0 =\emptyset$, nor the connectedness of the $D_i$'s. 
 Then we obtain~:
\begin{theorem}\label{gencas}~\\
 If $N=2$ and $\Omega$ is bounded, connected, satisfies the interior cone property and with a 
 piecewise $\mathcal C^{1,+}$ boundary, 
then the assertion of Theorem~\ref{Bernard}  holds with all 
 $a_i$'s strictly positive. \\
Moreover any minimizing partition $\mathcal D$ admits an open regular 
 connected representative.
\end{theorem}

{\bf Proof}~\\
Assume that some of the $a_i$ vanish or, in other words, that  ${\bf k}_0\neq\emptyset$. Let us denote
$$\mathcal W_0=\Omega\setminus \bigcup_{i\not\in{\bf k}_0}D_i$$
the nodal set. Then the 2--dimensional measure of $\mathcal W_0$ is
positive,
 since it contains the supports of $\widetilde\phi_i$ for all $i\in {\bf k}_0$.

Now, as already remarked, inequalities (I1) and (I2) are still
 available
 when we discard all the vanishing components and we take
 $U=(u_i)_{i\not\in{\bf k}_0}$.
 Hence Theorem~\ref{asym} applies and, as a direct consequence,
 we find that the zero set $\cal Z$, being the union of a finite number of
 closed $\mathcal C^{1,+}$
 curves has vanishing measure in $\mathbb R^2$, a contradiction.
This also implies that the partition associated with the $D_i$ is strong and
that the zero set is indeed the nodal set $N(\mathcal D)$ as defined in \eq{assclset}.

Now assume by contradiction some elements of the partition not be
connected. 
We can anyway choose a set of  first eigenfunctions having each a
connected support,
 but we clearly have an open, non empty set of multiplicity zero
 points. 
This contradicts Corollary \ref{no_vanishing}.\\
\begin{remark}~\\
Note that Theorem \ref{gencas} completes the proof 
 of Theorem \ref{thanyreg}. 
\end{remark}

{\bf Proof of Theorem \ref{partnod}}~\\
But we also  get the proof of Theorem \ref{partnod}
 in the following way. If the graph associated to $\widetilde D$  is bipartite
 we can find $\epsilon_i = \pm 1$ satisfying 
 $\epsilon_i\epsilon_j =-1$ if $D_i\sim D_j$  and such that $u:=
 \sum_i \epsilon_i a_i \widetilde \varphi_i$ is in $W_0^1(\Omega)$ 
 and satisfies 
\begin{equation}\label{valpr}
(-\Delta + V) u =\mathfrak L_k(\Omega) u
\end{equation}
 in $\Omega \setminus \mathcal Z_3$. But we have proven that  $\mathcal Z_3$
 consists of isolated  points (which cannot be the support of
 a distribution in $W^{-1}(\Omega)$). Hence  \eqref{valpr} is satisfied in $\Omega$
 and  $u$ is actually an eigenfunction of $H(\Omega)$
 corresponding to $\mathfrak L_k(\Omega)$. We refer to Appendix
 \ref{sregcase} for a complementary discussion.